\documentclass[preprint,3p,12pt]{elsarticle}
\usepackage{mathrsfs}
\usepackage{amsmath}
\usepackage{stmaryrd}
\usepackage{bbding}
\usepackage{dcolumn}
\usepackage{graphicx}
\usepackage{amsfonts}
\usepackage{amssymb}
\usepackage{psfrag}
\usepackage{wrapfig}
\usepackage{subfigure}
\usepackage{makeidx}
\usepackage{bm}
\usepackage{epsf}
\usepackage{color}
\usepackage{epsfig}
\usepackage{setspace}
\usepackage{epstopdf}
\usepackage{psfrag}
\usepackage{multirow}
\usepackage{diagbox}
\usepackage{verbatim}

\usepackage{makecell}
\usepackage{float}
\usepackage{esint}
\usepackage{comment}
\usepackage{movie15}
\usepackage{hyperref}
\usepackage[ruled,linesnumbered]{algorithm2e}
\usepackage{geometry}
\usepackage{appendix}
\usepackage{soul}

\newcommand{\mathsym}[1]{{}}
\newcommand{\unicode}[1]{{}}
\newcommand{\Rmnum}[1]{\uppercase\expandafter{\romannumeral #1}}
\newtheorem{remark}{Remark}[section]

\begin{document}
	\title{Treatment of Wall Boundary Conditions in High-Order Compact Gas-Kinetic Schemes}
	
	\author[HKUST1]{Jiawang Zhang}
	\ead{jzhangiw@connect.ust.hk}
	
	\author[XJTU]{Xing Ji}
	\ead{jixing@xjtu.edu.cn}
	
	\author[HKUST1,HKUST2,HKUST3]{Kun Xu\corref{cor1}}
	\ead{makxu@ust.hk}

	\address[HKUST1]{Department of Mathematics, Hong Kong University of Science and Technology, Clear Water Bay, Kowloon, Hong Kong}
	\address[XJTU]{State Key Laboratory for Strength and Vibration of Mechanical Structures, Shaanxi Key Laboratory of Environment and Control for Flight Vehicle, School of Aerospace Engineering, Xi'an Jiaotong University, Xi'an, China}
	\address[HKUST2]{Department of Mechanical and Aerospace Engineering, Hong Kong University of Science and Technology, Clear Water Bay, Kowloon, Hong Kong}
	\address[HKUST3]{Shenzhen Research Institute, Hong Kong University of Science and Technology, Shenzhen, China}
	\cortext[cor1]{Corresponding author}
	
	\begin{abstract}
		The boundary layer represents a fundamental structure in fluid dynamics, where accurate boundary discretization significantly enhances computational efficiency. This paper presents a third-order boundary discretization for compact gas-kinetic scheme (GKS). Wide stencils and curved boundaries pose challenges in the boundary treatment for high-order schemes, particularly for temporal accuracy. By utilizing a time-dependent gas distribution function, the GKS simultaneously evaluates fluxes and updates flow variables at cell interfaces, enabling the concurrent update of cell-averaged flow variables and their gradients within the  third-order compact scheme. The proposed one-sided discretization achieves third-order spatial accuracy on boundary cells by utilizing updated flow variables and gradients in the discretization for non-slip wall boundary conditions.  High-order temporal accuracy on boundary cells is achieved through the GKS time-dependent flux implementation with multi-stage multi-derivative methodology. Additionally, we develop exact no-penetration conditions for both adiabatic and isothermal wall boundaries, with extensions to curved mesh geometries to fully exploit the advantages of high-order schemes. Comparative analysis between the proposed one-sided third-order boundary scheme, third-order boundary scheme with ghost cells, and second-order boundary scheme demonstrates significant performance differences for the third-order compact GKS. Results indicate that lower-order boundary cell treatments yield substantially inferior results, while the proposed third-order treatment demonstrates superior performance, particularly on coarse grid configurations.

	\end{abstract}
	\maketitle
	\section{Introduction}
	
Boundary treatment in computational fluid dynamics (CFD) plays a crucial role in ensuring the numerical simulation accurately represents real-world flow physics. It provides the necessary conditions at the boundaries of the computational domain to properly simulate how the fluid interacts with surfaces and other domain limits. This treatment ensures that physical principles like non-slip conditions at walls, conservation of mass and energy, and proper flow entrance and exit conditions are correctly enforced. Proper boundary treatment is essential for maintaining numerical stability and solution accuracy, particularly in capturing important flow features like boundary layers and wake regions. Without appropriate boundary treatment, the simulation may produce unrealistic results, develop numerical instabilities, or fail to converge to a physical solution. This is especially important in cases involving complex geometries, high-order numerical schemes, or multi-physics problems where accurate boundary representation directly impacts the overall solution quality.
	
	The gas-kinetic scheme (GKS) is a computational fluid dynamics algorithm founded on kinetic equation \cite{BGK}. This method employs a time-dependent gas distribution function at cell interfaces to evaluate numerical fluxes \cite{xu1994, xu2001}. Through the implementation of multi-stage multi-derivative methodology, the scheme achieves fourth-order temporal accuracy using a two-stage process \cite{pan2016}. Beyond flux evaluation, the scheme uniquely enables the updating of flow variables at cell interfaces through the time-dependent gas distribution function, allowing simultaneous updates of cell-averaged flow variables and their gradients. This evolved variable approach facilitates the construction of compact high-order GKS on both structured and unstructured meshes \cite{ji2018,ji-hweno,zhao-acoustic}.
The compact GKS differs from other compact schemes such as the Lele-type compact scheme \cite{lele-compact, wang-cls}, Discontinuous Galerkin (DG) method \cite{dg-reed, dg-cockburn}, and correction procedure via reconstruction (CPR) method \cite{cpr, fr}. While Lele-type compact schemes feature implicit connections between flow variables and their derivatives within compact stencils, and DG methods employ implicit governing equations for multiple degrees of freedom (DOFs), these approaches utilize distinct evolution models for updating equivalent cell-averaged gradients. The versatility of GKS extends beyond Euler and Navier-Stokes solutions, finding applications in chemical reaction modelling  \cite{lian1}, shallow water equations \cite{zhao-swe}, magneto-hydrodynamics \cite{pu-plasma}, and turbulent flow simulations \cite{cao}.

The implementation of boundary conditions frequently relies on the ghost cell/node method, where the accurate determination of ghost cell states is crucial for solution quality and precision. Tan introduced the inverse Lax-Wendroff method to derive normal derivatives from tangential and temporal derivatives at boundaries \cite{invLW}, enabling precise ghost state implementation. Pointsot developed characteristic boundary conditions for the Navier-Stokes equations by solving characteristic systems at boundaries rather than constructing ghost cells \cite{nscbc}. Krivodonova proposed a high-order solid wall treatment for Euler flow by modifying flow patterns around solid walls to incorporate realistic geometry \cite{solidwall}, though this approach presents challenges when extended to viscous flow computations with complex geometries.

High-order wall boundary treatment necessitates high-order geometrical approximation of wall boundaries. Several high-order schemes incorporate high-order curved mesh computations, including high-order spectral methods \cite{sv-curve, sd-curve}, DG methods \cite{dg-curve, rdg-curve}, and compact GKS  \cite{ji-curve}. Research has demonstrated that inadequate boundary approximation compromises the accuracy and efficiency of compact high-order schemes \cite{bassi}. Studies also indicate that curved mesh representation of wall boundaries can reduce nonphysical entropy generation near walls  \cite{bn-effect}, while straight mesh lines can lead to decreased lift and increased drag due to artificial surface roughness \cite{lift-curve}.

This research focuses on developing an efficient, high-order wall boundary treatment by incorporating physical properties of non-slip wall boundaries—specifically zero velocity and zero normal pressure gradient—combined with updated flow variables and their gradients from interior cells. These elements are integrated to construct a one-sided third-order compact scheme for boundary cells. High-order temporal accuracy is achieved through the GKS time-dependent flux using the multi-stage multi-derivative method \cite{s2o4}. The constrained least-squares method enhances the stability of the one-sided reconstruction scheme \cite{liwanai}, with this high-order wall boundary treatment extended to curved meshes.

The paper is structured as follows: Section 2 presents the basic framework of the compact third-order GKS, Section 3 details the boundary treatment methodology, Section 4 demonstrates numerical examples with adiabatic and isothermal walls, and the final section provides concluding remarks.

	\section{Compact Finite Volume Gas-kinetic Scheme}
	Gas-kinetic solver is constructed from the gas-kinetic BGK equation,
	\begin{equation}\label{bgk}
		f_t+\mathbf{u}\cdot \nabla f = \frac{g-f}{\tau},
	\end{equation}
	where $f = f(\mathbf{x}, t, \mathbf{u}, \xi)$ is the gas distribution function, which is a function of space $\mathbf{x}$, time $t$, particle velocity $\mathbf{u}$ and internal  variable $\xi$. $\tau$ is the relaxation time from $f$ to its equilibrium state $g$. The equilibrium state $g$ is a Maxwellian distribution,
	\begin{equation}
		g = \rho (\frac{\lambda}{\pi})^{\frac{K + 2}{2}}e^{-\lambda [(\mathbf{u} - \mathbf{U})^2 + \xi^2]}. \nonumber
	\end{equation}
	
	Due to the conservation of mass, momentum and energy during particle collision process, the collision term in Eq. \ref{bgk} should satisfy the compatibility condition,
	\begin{equation} \label{comp_cond}
		\int \frac{g-f}{\tau} \psi \mathrm{d} \Xi = 0,
	\end{equation}
	where $\psi = (1, \mathbf{u}, \frac{1}{2}(\mathbf{u}^2 + \xi^2))^{T}$, $\mathrm{d} \Xi = \mathrm{d}u\mathrm{d}v\mathrm{d}\xi_1 \cdots \mathrm{d}\xi_K$, $K$ is the number of internal degree of freedom, $K=3$ for diatomic gas in 2D case.
	The conservative variables and fluxes can been obtained from the gas distribution function,
	\begin{equation} \label{consvar}
		\mathbf{W} = \int f \psi \mathrm{d} \Xi,
	\end{equation}
	\begin{equation} \label{flux}
		\mathbf{F} = \int f \psi \mathbf{u} \mathrm{d} \Xi,
	\end{equation}
	where $\mathbf{W}=\{\rho, \rho U, \rho V, \rho E\}^{T}$ is the conservative variables and $\mathbf{F}$ is the corresponding flux.

	\subsection{Finite volume formulation}
	Take the moments of the BGK model Eq. \ref{bgk} in phase space and integrate it in a finite volume $\Omega_i$,
	\begin{equation}
		\int_{\Omega_i}\int (f_t+\mathbf{u}\cdot \nabla f)\psi \mathrm{d}\Xi\mathrm{d} V = \int_{\Omega_i}\int \frac{g-f}{\tau}\psi \mathrm{d}\Xi\mathrm{d} V, \nonumber
	\end{equation}
	due to the relations Eq. \ref{comp_cond}, Eq. \ref{consvar} and Eq. \ref{flux}, the integral form is obtained
	\begin{equation} \label{int_form}
		\int_{\Omega_i} \mathbf{W}_t \mathrm{d} V + \int_{\Omega_i} \nabla \cdot \mathbf{F} \mathrm{d} V = 0, \nonumber
	\end{equation}
	The above integral form is discretized by finite volume method (FVM),
	\begin{equation} \label{semi_discrete}
		\frac{\mathrm{d} \overline{\mathbf{W}}_i}{\mathrm{d} t} = -\frac{1}{|\Omega_i|}\int_{\Omega_i} \nabla \cdot \mathbf{F} \mathrm{d} V.
	\end{equation}
	From Gaussian theorem, the semi-discrete form Eq. \ref{semi_discrete} is written as,
	\begin{equation} \label{res}
		\frac{\mathrm{d} \overline{\mathbf{W}}_i}{\mathrm{d} t}
		= \mathcal{L}(\mathbf{W}_i)
		=-\frac{1}{|\Omega_i|} \int_{\partial \Omega_i} \mathbf{F} \cdot \mathbf{n} \mathrm{d}s
		= -\frac{1}{|\Omega_i|} \sum_{p=1}^{N_f}\int_{\Gamma_{ip}} \mathbf{F} \cdot \mathbf{n}_p \mathrm{d}s,
	\end{equation}
	where $|\Omega_i|$ is the volume of the control volume, $\partial \Omega_i$ is the boundary of the control volume, which is expressed as,
	\begin{equation}
		\partial \Omega_i = \bigcup_{p=1}^{N_f}\Gamma_{ip}, \nonumber
	\end{equation}
	where $\Gamma_{ip}$ is the neighboring interface of the cell $\Omega_i$, $N_f$ is the number of cell interface.
	Numerical quadrature is used to evaluate the surface integral of fluxes,
	\begin{equation}
		\int_{\Gamma_{ip}} \mathbf{F} \cdot \mathbf{n}_p \mathrm{d}s = |\Gamma_{ip}|\sum_{k=1}^{N_g} \omega_k \mathbf{F}(\mathbf{x}_{p,k},t)\cdot \mathbf{n}_{p,k}, \nonumber
	\end{equation}
	where $|\Gamma_{ip}|$ is the area of the mesh face, $\mathbf{x}_{p,k}$ is the Gaussian quadrature points and $\omega_k$ is the Gaussian weight, $N_g$ is the number of Gaussian points.
	\subsection{Gas-kinetic solver}
	Gas-kinetic scheme is based on the exact integral solution of the BGK equation Eq. \ref{bgk},
	\begin{equation} \label{integral_sol}
		f(\mathbf{x}, t, \mathbf{u}, \xi) = \frac{1}{\tau}\int_{0}^{t} g(\mathbf{x}', t', \mathbf{u}, \xi)e^{-\frac{t-t'}{\tau_n}}\mathrm{d}t'+e^{\frac{-t}{\tau_n}}f_0(\mathbf{x}_0, \mathbf{u}),
	\end{equation}
	where $\mathbf{x} = \mathbf{x}' + \mathbf{u} (t-t')$ is the trajectory of the particle motion, $\mathbf{x} = \mathbf{0}$ is the location of cell interface,
$f_0$ is the gas distribution function at the beginning of each time step ($t=0$), $\mathbf{x}_0$ is the initial position of the concerned particle.
For the gas-kinetic solver, the equilibrium term $g(\mathbf{x}', t', \mathbf{u}, \xi)$ and initial distribution function $f_0(\mathbf{x}_0, \mathbf{u})$ need to be constructed.
	To construct a time evolution solution of gas distribution function at a cell interface, the following notations are introduced first,
	\begin{equation}
		\mathbf{a} = \frac{1}{g}\frac{\partial g}{\partial \mathbf{x}}, A = \frac{1}{g}\frac{\partial g}{\partial t}, \nonumber
	\end{equation}
	where $\mathbf{a}$ and $A$ are the spatial and temporal derivative of the equilibrium distribution function.
	According to Taylor expansion, the spatial and temporal distribution of the equilibrium state is
	\begin{equation} \label{g}
		\begin{aligned}
			g(\mathbf{x}', t', \mathbf{u}, \xi) = g(\mathbf{0}, 0, \mathbf{u}, \xi)(1-\mathbf{a} \cdot \mathbf{u}(t-t') + At').
		\end{aligned}
	\end{equation}
	The initially spatial distribution function around the interface is
	\begin{equation}\label{f0}
		\begin{aligned}
			f_0(\mathbf{x}_0, 0, \mathbf{u}, \xi)
			&=g^{l}(\mathbf{0}, 0, \mathbf{u}, \xi)(1-(t+\tau)\mathbf{a}^{l/r} \cdot \mathbf{u} - \tau A^{l/r})H(u) \\
			&+g^{r}(\mathbf{0}, 0, \mathbf{u}, \xi)(1-(t+\tau)\mathbf{a}^{l/r} \cdot \mathbf{u} - \tau A^{l/r})(1-H(u)).
		\end{aligned}
	\end{equation}
	Substituting Eq. \ref{g} and Eq. \ref{f0} into the integral solution Eq. \ref{integral_sol}, the time evolved distribution function at the interface can be obtained
	\begin{equation} \label{f_inter}
		\begin{aligned}
			f(\mathbf{0}, t, \mathbf{u}, \xi) &= C_1g^c+C_2\mathbf{a}^c\cdot \mathbf{u} g^c + C_3A^cg^c \\
			&+ C_4[g^lH(u) + g^r(1-H(u))] \\
			&+ C_5[g^l\mathbf{a}^l\cdot \mathbf{u}H(u) + g^r\mathbf{a}^r\cdot \mathbf{u}(1-H(u))] \\
			&+ C_6[g^lA^lH(u) + g^rA^r(1-H(u))],
		\end{aligned}
	\end{equation}
	with the coefficients $C_1 = 1-e^{-t/\tau_n}$, $C_2 = (t+\tau)e^{-t/\tau_n}-\tau$, $C_3 = t-\tau+\tau e^{-t/\tau_n}$, $C_4 = e^{-t/\tau_n}$, $C_5 = -e^{-t/\tau_n}(\tau + t)$, $C_6 = -\tau e^{-t/\tau_n}$.
	For smooth flow, the time-dependent solution Eq. \ref{f_inter} can be simplified as
	\begin{equation} \label{gks-smooth}
		f(\mathbf{0}, t, \mathbf{u}, \xi) = g^c - \tau (\mathbf{a}^c \cdot \mathbf{u} + A^c)g^c + A^c g^ct.
	\end{equation}
	For viscous flow, the physical collision time is determined by
	\begin{equation}
		\tau = \frac{\mu}{p},  \nonumber
	\end{equation}
	where $\mu$ is the dynamic viscosity. To properly capture the discontinuity with additional numerical dissipation, the numerical collision time is modified as
	\begin{equation}
		\tau_n = \frac{\mu}{p} + C\frac{|p_l - p_r|}{|p_l + p_r|}\Delta t,  \nonumber
	\end{equation}
	where $\Delta t$ is the time step, $p_l$ and $p_r$ are the pressure on the left and right side of the interface.
	\begin{remark}
		Unit Prandtl number $Pr = 1$ is associated with the single relaxation BGK model. To simulate fluid with non-unit Prandtl number, e.g. $Pr \approx 0.72$ for air, the energy flux is modified according to the real Prandtl number \cite{xu2001},
		\begin{equation}
			F_{\rho E}^{new} = F_{\rho E} + (\frac{1}{Pr} - 1)q, \nonumber
		\end{equation}
		where $q$ is the heat flux calculated by
		\begin{equation}
			q = \int \frac{1}{2} (u-U)[(\mathbf{u}-\mathbf{U})^2 + \xi^2]f\mathrm{d}\Xi.  \nonumber
		\end{equation}
	\end{remark}
	
	\subsection{Gas-kinetic solver for isothermal boundary} \label{gks_iso}
	To deal with isothermal boundary condition, a kinetic boundary condition was constructed by Li \cite{kineticbn}, where a Maxwellian distribution was imposed on the outer side of the wall boundary. Here a time-dependent distribution function is imposed on the outer side. Based on the reconstructed left state at the inner side of the boundary face, the inner-side distribution function is
	\begin{equation}
		f^L = g^L(1-(t+\tau)\mathbf{a}^L \cdot \mathbf{u} - \tau A^L).
	\end{equation}
	Supposing the outer-side distribution function is
	\begin{equation}  \label{iso_fR}
		f^R = g^R(1 + A^Rt),
	\end{equation}
	where
	\begin{equation}
		\int g^R A^R \psi \mathbf{d} \Xi = \frac{\partial \mathbf{W}^R}{\partial t}, \nonumber
	\end{equation}
	$A^R$ is obtained from the Taylor expansion of a Maxwellian and has the form
	\begin{equation}
		A^R = A_1^R + A_2^Ru + A_3^Rv + A_4^R \frac{1}{2} (u^2 + v^2 + \xi^2). \nonumber
	\end{equation}
	From Ref. \cite{xu2001}, $A^R$ is determined by
	\begin{equation}
		\begin{aligned}
			A_4^R &= \frac{8\lambda^2}{K + 2}(\frac{\partial E}{\partial t} - U\frac{\partial U}{\partial t} -V \frac{\partial V}{\partial t}), \\
			A_3^R &= 2\lambda(\frac{\partial V}{\partial t} - \frac{V}{2\lambda}A_4^R), \\
			A_2^R &= 2\lambda(\frac{\partial U}{\partial t} - \frac{U}{2\lambda}A_4^R), \\
			A_1^R &= \frac{1}{\rho}\frac{\partial \rho}{\partial t} - U A_2^R - V A_3^R - E A_4^R.
		\end{aligned} \nonumber
	\end{equation}
	
	According to the conditions $U = 0, V = 0, E = C_vT + 0.5 (U^2 + V^2) = C_v T_w$ at the non-slip isothermal wall, we have
	\begin{equation}
		\begin{gathered}
			A_1^R = \frac{1}{\rho}\frac{\partial \rho}{\partial t}, \\
			A_{2-4}^R = 0.
		\end{gathered} \nonumber
	\end{equation}
	
	There are two unknowns in Eq. \ref{iso_fR}, $\rho^R$ and $A_1^R$, which are determined by mass no-penetration condition at $t=0$ and $t=\Delta t$,
	\begin{equation}
		\int u(f^L(t = 0)H(u) + f^R(t = 0)(1-H(u)))\mathbf{d} \Xi = 0, \nonumber
	\end{equation}
	\begin{equation}
		\int u(f^L(t = \Delta t)H(u) + f^R(t = \Delta t)(1-H(u)))\mathbf{d} \Xi = 0, \nonumber
	\end{equation}
	leading to
	\begin{equation}
		\int_{u<0} ug^R \mathbf{d} \Xi = - \int_{u>0}ug^L(1-\tau (\mathbf{a} \cdot \mathbf{u} + A^L)) \mathbf{d} \Xi, \nonumber
	\end{equation}
	\begin{equation}
		\int_{u<0} ug^R A^R\mathbf{d} \Xi = \int_{u>0}ug^L(\mathbf{a} \cdot \mathbf{u}) \mathbf{d} \Xi. \nonumber
	\end{equation}
	
	Then the unknowns $\rho^R$ and $A^R$ can be determined. And the flux across the isothermal interface is calculated by
	\begin{equation} \label{iso_flux}
		\mathbf{F} = \int_{u>0} u f^L \psi \mathbf{d} \Xi + \int_{u<0} u f^R \psi \mathbf{d} \Xi.
	\end{equation}
	
	The flux evaluated by Eq. \ref{iso_flux} can guarantee mass no-penetration rigidly, which have been verified by the numerical cases in Section 4.
	
	\subsection{Direct evolution of the cell-averaged gradients of flow variables}
	Thanks to the time-dependent distribution function in  GKS Eq. \ref{f_inter}, the conservative variables at the cell interface can also be evaluated by
	\begin{equation}
		\mathbf{W}(\mathbf{x}_{p,k}, \Delta t) = \int f(\mathbf{x}_{p,k}, \Delta t, \mathbf{u}, \xi)\psi \mathrm{d}\Xi. \nonumber
	\end{equation}
	From Gaussian theorem, the cell-averaged gradients of the flow variables are obtained from the updated conservative variables at the Gaussian points,
	\begin{equation}
		\overline{\nabla \mathbf{W}}^{n+1} = \frac{1}{|\Omega_i|}\int_{\Omega_i} \nabla \mathbf{W}(\Delta t) \mathrm{d} V = \frac{1}{|\Omega_i|}\sum_{p=1}^{N_f} |\Gamma_{ip}|\sum_{k=1}^{N_g}\omega_k \mathbf{W}(\mathbf{x}_{p,k}, \Delta t) \mathbf{n}_{p,k}.
	\end{equation}
	\subsection{Two-stage fourth-order temporal scheme}
	The two-stage fourth-order temporal scheme is adopted to discretize the semi-discrete form Eq. \ref{res},
	\begin{equation}
		\begin{gathered}
			\mathbf{W}_i^{n+1/2} = \mathbf{W}_i^{n} + \frac{1}{2}\Delta t \mathcal{L}(\mathbf{W}_i^{n}) + \frac{1}{8}\Delta t^2 \mathcal{L}^{(1)}(\mathbf{W}_i^n), \\
			\mathbf{W}_i^{n+1} =\mathbf{W}_i^{n} + \Delta t \mathcal{L}(\mathbf{W}_i^{n}) + \frac{1}{6}\Delta t^2 (\mathcal{L}^{(1)}(\mathbf{W}_i^n) + 2\mathcal{L}^{(1)}(\mathbf{W}_i^{n+1/2})),
		\end{gathered}
	\end{equation}
	where $\mathcal{L}(\mathbf{W}_i^{n})$ and $\mathcal{L}^{(1)}(\mathbf{W}_i^n)$ are net flux and its first-order temporal derivative of cell $\Omega_i$ at time $t=t_n$, which are evaluated by numerical integral,
	\begin{equation}
		\mathcal{L}(\mathbf{W}_i^{n}) = \frac{1}{|\Omega_i|} \sum_{p=1}^{N_f}|\Gamma_{ip}|\sum_{k=1}^{N_g} \omega_k \mathbf{F}(\mathbf{x}_{p,k},t_n)\cdot \mathbf{n}_{p,k}, \nonumber
	\end{equation}
	\begin{equation}
		\mathcal{L}^{(1)}(\mathbf{W}_i^{n}) = \frac{1}{|\Omega_i|} \sum_{p=1}^{N_f}|\Gamma_{ip}|\sum_{k=1}^{N_g} \omega_k \partial_t \mathbf{F}(\mathbf{x}_{p,k},t_n)\cdot \mathbf{n}_{p,k}, \nonumber
	\end{equation}
	where $\mathbf{F}(\mathbf{x}_{p,k},t_n)$ is the flux through the Gaussian point $\mathbf{x}_{p,k}$ at time $t=t_n$ and $\partial_t \mathbf{F}(\mathbf{x}_{p,k},t_n)$ is the temporal derivative of the flux. To evaluate the flux and its temporal derivative, the following notation is introduced,
	\begin{equation} \label{flux_dt}
		\mathbb{F}(\mathbf{x}_{p,k}, \Delta t) = \int_{0}^{\Delta t}\mathbf{F}(\mathbf{x}_{p,k}, t) \mathrm{d}t.
	\end{equation}
	The Taylor expansion of the flux at the interface with second order truncation error is
	\begin{equation} \label{taylor_flux}
		\mathbf{F}(\mathbf{x}_{p,k}, t) = \mathbf{F}(\mathbf{x}_{p,k}, 0) + t\partial_t \mathbf{F}(\mathbf{x}_{p,k}, t).
	\end{equation}
	Substituting the Taylor expansion Eq. \ref{taylor_flux} into the integral Eq. \ref{flux_dt},
	\begin{equation}
		\mathbf{F}(\mathbf{x}_{p,k}, 0)\Delta t + \frac{1}{2} \partial_t \mathbf{F}(\mathbf{x}_{p,k}, 0)\Delta t^2 = \mathbb{F}(\mathbf{x}_{p,k}, \Delta t),  \nonumber
	\end{equation}
	\begin{equation}
		\frac{1}{2}\mathbf{F}(\mathbf{x}_{p,k}, 0)\Delta t + \frac{1}{8} \partial_t \mathbf{F}(\mathbf{x}_{p,k}, 0)\Delta t^2 = \mathbb{F}(\mathbf{x}_{p,k}, \Delta t / 2).  \nonumber
	\end{equation}
	the flux and its temporal derivative at the interface at the beginning of each time step can be derived,
	\begin{equation}
		\mathbf{F}(\mathbf{x}_{p,k}, 0) = (4\mathbb{F}(\mathbf{x}_{p,k}, \Delta t / 2) - \mathbb{F}(\mathbf{x}_{p,k}, \Delta t)) / \Delta t, \nonumber
	\end{equation}
	\begin{equation}
		\partial_t \mathbf{F}(\mathbf{x}_{p,k}, 0) = 4(\mathbb{F}(\mathbf{x}_{p,k}, \Delta t) - 2 \mathbb{F}(\mathbf{x}_{p,k}, \Delta t / 2)) / \Delta t^2. \nonumber
	\end{equation}
	The time step is determined by
	\begin{equation}
		\Delta t = \mathrm{CFL} \times \min \{ \frac{\Delta h}{|\mathbf{u}| + c}, \frac{\Delta h^2}{4 v}\}, \nonumber
	\end{equation}
	where $\rm{CFL}$ is the Courant-Friedrichs-Lewy number, $\Delta h$ is the characteristic length of the cell, $\mathbf{u}$ is the velocity, $c$ is the sound speed, $v=\mu / \rho$ is the kinematic viscosity.
	\section{Spatial Reconstruction}
	From the evolved cell-averaged flow variables and their spatial gradients by the CGKS framework, the Hermite weighted essentially non-oscillatory (HWENO) \cite{Qiu-hweno} can be applied to construct the third-order compact scheme. For wall boundary cells, the physical properties of wall boundary are directly used in the reconstruction instead of constructing ghost cells.
	\subsection{HWENO Reconstruction}  \label{hweno}
	The HWENO reconstruction has been applied to construct compact high-order gas-kinetic scheme in Ref. \cite{ji-hweno}, where a quadratic polynomial is constructed by Hermite polynomial in a compact stencil consisting of only the face-neighboring cells, and a linear polynomial is constructed by Green Gauss method. In order to capture discontinuity, the linear and quadratic polynomial is weighted based on the discontinuity feedback (DF) factor proposed in Ref. \cite{ji-cf}. \\
	
	\subsubsection{Construction of quadratic polynomial}
\begin{figure}[!htbp]
		\centering
		\includegraphics[width = 0.3\columnwidth]{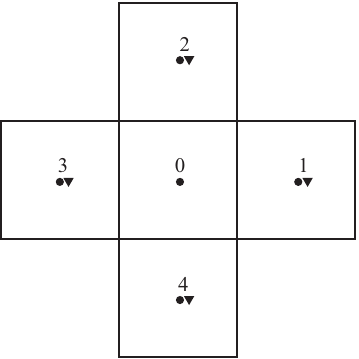}
		\caption{Reconstruction stencil for compact third order scheme}
		\label{fig:stencil}
	\end{figure}

	To achieve third-order accuracy, a quadratic polynomial is constructed
	\begin{equation} \label{p_2}
		p^2(x,y) = \overline{Q}_{0} + \sum_{k=1}^{5}a_k\phi_k(x,y),
	\end{equation}
	where $\overline{Q}_{0}$ is the cell-averaged value of $Q(x,y)$ over cell $\Omega$, $a_k$ are the unknown coefficients of the quadratic polynomial and $\phi_k(x,y), k=1,\cdots, 5$ are zero-mean basis functions,
	\begin{equation}
		\begin{aligned}
			\phi_1(x,y) &= x - x_0, \\
			\phi_2(x,y) &= y - y_0, \\
			\phi_3(x,y) &= (x-x_0)^2 - \frac{1}{|\Omega|}\iint_{\Omega} (x-x_0)^2\mathrm{d}x\mathrm{d}y, \\
			\phi_4(x,y) &= (x-x_0)(y-y_0) - \frac{1}{|\Omega|}\iint_{\Omega} (x-x_0)(y-y_0)\mathrm{d}x\mathrm{d}y, \\
			\phi_5(x,y) &= (y-y_0)^2 - \frac{1}{|\Omega|}\iint_{\Omega} (y-y_0)^2\mathrm{d}x\mathrm{d}y,
		\end{aligned} \nonumber
	\end{equation}
	where $(x_0,y_0)$ is the center of the concerned cell. The reconstruction stencil for quadrilateral cell is demonstrated in Figure \ref{fig:stencil}.
	
	To determine the unknown coefficients in Eq. \ref{p_2}, the following conditions should be satisfied,
	\begin{equation} \label{recon_cond}
		\begin{aligned}
			&\frac{1}{|\Omega_j|}\iint_{\Omega_j} p^2(x,y)\mathrm{d}x\mathrm{d}y = \overline{Q}_j, \\
			&\frac{1}{|\Omega_j|}\iint_{\Omega_j} \frac{\partial}{\partial x} p^2(x,y)\mathrm{d}x\mathrm{d}y = \overline{Q_x}_{,j}, \\
			&\frac{1}{|\Omega_j|}\iint_{\Omega_j} \frac{\partial}{\partial y} p^2(x,y)\mathrm{d}x\mathrm{d}y = \overline{Q_y}_{,j}, \\
		\end{aligned}
	\end{equation}
	where $j = 1,2,3,4$ indicates the face-neighboring cell, $\overline{Q}_j$, $\overline{Q_x}_{,j}$, $\overline{Q_y}_{,j}$ are the cell-averaged value, $x$-derivative and $y$-derivative of variable $Q$ over the cell $j$ respectively.  There are 12 constrains by Eq. \ref{recon_cond}, but only five unknowns in the quadratic polynomial Eq. \ref{p_2}. To solve this over-determined system and enhance the linear stability, the constrained least squares technique in Ref. \cite{liwanai} is adopted, where the cell-averaged values are constrained.
	
	\subsubsection{Reconstruction of linear polynomial}
	
	The linear polynomial is in the form
	\begin{equation} \label{p_1}
		p^1(x,y) = \overline{Q}_{0} + a_1 (x-x_0) + a_2(y-y_0).
	\end{equation}
	The unknown coefficients $a_1$ and $a_2$ are determined by divergence theorem,
	\begin{equation}
		\iint_{\Omega} \nabla p^1(x,y) \mathrm{d}x\mathrm{d}y = \oiint_{\partial \Omega} Q(x,y) \mathbf{n}\mathrm{d}s = \sum_{j=1}^{4} \iint_{\Gamma_j} Q(x,y) \mathbf{n}\mathrm{d}s \approx \sum_{j=1}^{4} |\Gamma_j|\overline{Q}_{0,j} \mathbf{n}_j, \nonumber
	\end{equation}
	where $|\Gamma_j|$ is the area of $j-th$ cell interface, $\mathbf{n}_j$ is the outer unit normal vector of $j-th$ interface, $\overline{Q}_{0,j}$ is the averaged value of variable $Q$ over the interface, which is approximated by
	\begin{equation}
		\overline{Q}_{0,j} = \frac{1}{2} (\overline{Q}_{0} + \overline{Q}_{j}). \nonumber
	\end{equation}
	
	To properly capture discontinuities, the discontinuity feedback factor is used to limit the coefficients in Eq. \ref{p_1}, which is defined as
	\begin{equation}
		\alpha_c = \prod_{j=1}^{N_f}\prod_{k=1}^{N_g} \alpha_{j,k}, \alpha_c \in (0,1], \nonumber
	\end{equation}
	where $\alpha_{j,k}$ is the DF factor at $k-th$ Gaussian point of $j-th$ interface, as in Ref. \cite{sliding}, DF factor is evaluated by
	\begin{equation}
		\begin{aligned}
			&\alpha_{j,k} = \frac{1}{1+A^2}, \\
			&A = \frac{|p^l - p^r|}{p^l} + \frac{|p^l - p^r|}{p^r} + (Ma_n^l - Ma_n^r)^2 + (Ma_t^l - Ma_t^r)^2,
		\end{aligned} \nonumber
	\end{equation}
	where $p$ is the pressure, $Ma_n$ and $Ma_t$ are the normal and tangential Mach number at the Gaussian point. Then the linear polynomial Eq. \ref{p_1} is limited by
	\begin{equation}
		\widetilde{p}^1(x,y) = \overline{Q}_{0} + \alpha_c \cdot (a_1 (x-x_0) + a_2(y-y_0)).
	\end{equation}
	
	\subsubsection{Non-linear weighting of quadratic and linear polynomial}
	
	To obtain third-order accuracy in smooth region and improve the robustness in discontinuous region, we derive an adaptive selection between the quadratic polynomial $p^2$ and the limited nonlinear polynomial $\widetilde{p}^1$ by DF factor,
	\begin{equation}  \label{nonlinear-p2}
		P^2 = H[\alpha_c - \Theta]p^2 + (1-H[\alpha_c - \Theta])\widetilde{p}^1,
	\end{equation}
	where $H[x]$ is Heavisid step function, $\Theta$ is a threshold for the selection of quadratic polynomial, which takes $\Theta = 0.9$ in this work.
	\subsubsection{Reconstruction on reference coordinate system}
	When simulating viscous flow with high Reynolds number, the cells near wall boundary may have large aspect ratio, leading to a large condition number of the constrained least-squares method. To enhance the stability of reconstruction on cells with large aspect ratio, the reconstruction is conducted on reference coordinate system as in Ref \cite{zhao-acoustic}. For a quadrilateral cell, the transformation from physical coordinate system $x-y$ to reference coordinate system $\xi - \eta$ is,
	\begin{figure}[!htbp]
		\centering
		\includegraphics[width = 0.75\columnwidth]{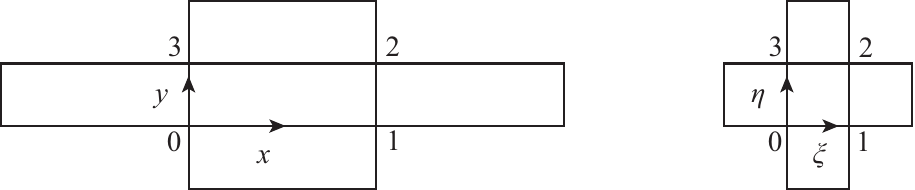}
		\caption{Transformation from physical system to reference coordinate system}
		\label{fig:ref_sys}
	\end{figure}
	\begin{equation}  \label{jacobi}
		\begin{pmatrix}
			x \\
			y
		\end{pmatrix} =
		\begin{pmatrix}
			x_0 \\
			y_0
		\end{pmatrix} +
		\mathbf{J}
		\begin{pmatrix}
			\xi \\
			\eta
		\end{pmatrix},
	\end{equation}
	where $\mathbf{J}$ is the Jacobian matrix, $(x_j, y_j), j = 0, 1, 2, 3$ are the coordinates of the vertexes of the quadrilateral. With this transformation, the reconstruction stencil in physical system $S$ can be transformed to the reference system $\widetilde{S}$, the Jacobian matrix takes,
	\begin{equation}
		\mathbf{J} =
		\begin{pmatrix}
			x_1 - x_0 & x_3 - x_0 \\
			y_1 - y_0 & y_3 - y_0 \\
		\end{pmatrix}. \nonumber
	\end{equation}
	Flow variables and their derivatives in reference coordinate system are obtained from
	\begin{equation}
		\widetilde{Q} = Q, \nonumber
	\end{equation}
	\begin{equation}
		\begin{pmatrix}
			\widetilde{Q}_{\xi} \\
			\widetilde{Q}_{\eta}
		\end{pmatrix} = \mathbf{J}^{\mathrm{T}}
		\begin{pmatrix}
			Q_x \\
			Q_y
		\end{pmatrix}. \nonumber
	\end{equation}
	Then the reconstruction in Eq. \ref{recon_cond} can be conducted in reference coordinate system. The flow variables at Gaussian points can be interpolated in reference system after the reconstruction, then the variables in reference system should be transformed to physical system by
	\begin{equation}
		Q = \widetilde{Q}, \nonumber
	\end{equation}
	\begin{equation}
		\begin{pmatrix}
			Q_x \\
			Q_y
		\end{pmatrix} = \mathbf{J}^{\mathrm{-T}}
		\begin{pmatrix}
			\widetilde{Q}_{\xi} \\
			\widetilde{Q}_{\eta}
		\end{pmatrix}. \nonumber
	\end{equation}
	
	\subsection{Reconstruction of boundary cells} \label{bnRecon}
	The physical constrains of wall boundary are considered in the reconstruction process without constructing ghost cells, e.g., the local velocity $U=0, V=0$, and the normal gradient of pressure $\frac{\partial p}{\partial \mathbf{n}} = \mathbf{0}$ at non-slip walls, $T = T_w$ for isothermal wall. The reconstruction is based on conservative variables in this work, so the constraints of primitive variable on wall boundary should be transferred to those of conservative variables.
	\subsubsection{Non-slip adiabatic wall}
	Due to $U, V = 0$ at non-slip walls, the normal gradient of total energy at non-slip wall boundary is determined by that of pressure,
	\begin{equation}
		\frac{\partial \rho E}{\partial \mathbf{n}} = \frac{1}{\gamma - 1}\frac{\partial p}{\partial \mathbf{n}}
		+ \frac{\partial}{\partial \mathbf{n}} (\frac{1}{2} \rho (U^2 + V^2)) = \frac{1}{\gamma - 1}\frac{\partial p}{\partial \mathbf{n}}. \nonumber
	\end{equation}
	The normal gradient of pressure can be expressed by density and temperature,
	\begin{equation}
		\frac{\partial p}{\partial \mathbf{n}} = \frac{\partial}{\partial \mathbf{n}} (\rho R_g T)
		= R_g (\frac{\partial \rho}{\partial \mathbf{n}} T + \rho \frac{\partial T}{\partial \mathbf{n}}). \nonumber
	\end{equation}
	For non-slip walls, the normal gradient of pressure is zero in boundary layer,
	\begin{equation}
		\frac{\partial \rho}{\partial \mathbf{n}} T + \rho \frac{\partial T}{\partial \mathbf{n}} = 0. \nonumber
	\end{equation}
	For adiabatic wall, the normal gradient of temperature is zero $\frac{\partial T}{\partial \mathbf{n}} = 0$, which lead to
	\begin{equation} \label{rho-adia}
		\frac{\partial \rho}{\partial \mathbf{n}} = 0.
	\end{equation}
	Above all, the constrained conditions for conservative variables of non-slip adiabatic wall read
	\begin{equation} \label{adia-cond}
		\begin{aligned}
			\frac{\partial \rho}{\partial \mathbf{n}} = 0,
			\rho U = 0,
			\rho V = 0,
			\frac{\partial \rho E}{\partial \mathbf{n}} = 0.
		\end{aligned}
	\end{equation}
	\subsubsection{Non-slip isothermal wall} \label{recon_iso}
	For non-slip isothermal wall, the velocity and normal gradient of pressure are still zero at isothermal wall boundary, so the constrains for $\rho U, \rho V, \rho E$ still hold,
	\begin{equation}
		\begin{aligned}
			\rho U = 0,
			\rho V = 0,
			\frac{\partial \rho E}{\partial \mathbf{n}} = 0.
		\end{aligned} \nonumber
	\end{equation}
	
	But the normal gradient of temperature is not zero and depends on the heat flux across isothermal boundary, so the constraint of density Eq. \ref{rho-adia} does not hold. To determine the density at isothermal wall, the momentum $\rho U, \rho V$ and total energy $\rho E$ are reconstructed firstly, once $\rho E$ at boundary is determined, the pressure at boundary can be obtained by
	\begin{equation}
		p_w = (\gamma - 1)(\rho E)_w. \nonumber
	\end{equation}
	From the equation of state, the density at boundary is determined by
	\begin{equation}
		\rho_w = \frac{p_w}{R_g T_w}. \nonumber
	\end{equation}
	Finally, density will be reconstructed at the boundary cells near isothermal wall. Then, the inner states at the left side of the isothermal boundary are all determined, and the isothermal gas-kinetic solver described in Section \ref{gks_iso} is used to evaluate the flux.
	
	\subsection{Curved mesh}
\begin{figure}[!htbp]
		\centering
		\subfigure[Green Gauss]{\includegraphics[width = 0.265\columnwidth]{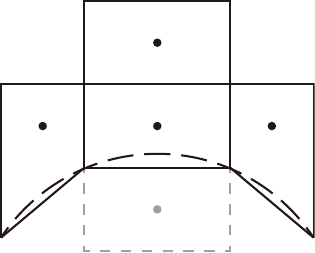}} \hspace{1pt}
		\subfigure[Third-order ghost cell]{\includegraphics[width = 0.265\columnwidth]{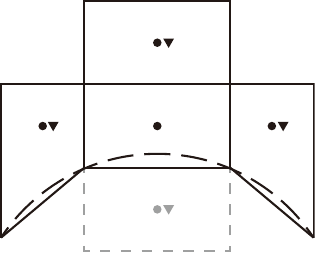}} \hspace{1pt}
		\subfigure[Third-order one-side]{\includegraphics[width = 0.265\columnwidth]{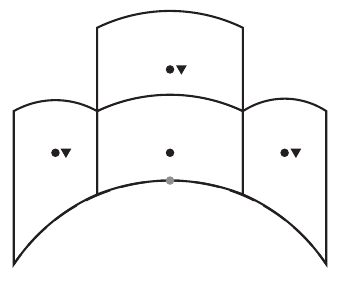}} \hspace{1pt}
		\subfigure{\includegraphics[width = 0.15\columnwidth, trim = 0 30 0 0, clip]{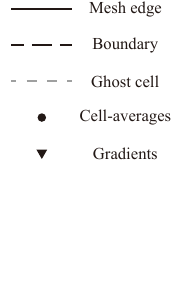}}
		\caption{Reconstruction stencil for (a) Green Gauss scheme, (b) compact third order scheme with ghost cells and (c) compact third order scheme with one-sided stencil}
		\label{fig:bn-sten}
	\end{figure}

	In this work, three types of boundary treatment are compared, and the reconstruction stencils are shown in Figure \ref{fig:bn-sten}. In the first method, the Green Gauss method is used for the reconstruction on wall boundary cells. In the second method, the third-order compact scheme with ghost cells is used, and the setting of the ghost state in the ghost cells is introduced in detail in Ref. \cite{ji-curve}. In the third method, the proposed boundary treatment introduced in Section \ref{bnRecon} is used with quadratic cells to match the quadratic reconstruction polynomial.
	
\begin{figure}[!htbp]
		\centering
		\includegraphics[width = 0.75\columnwidth]{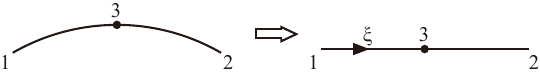}
		\caption{Isoparametric transformation for quadratic curve}
		\label{fig:quad-curve}
	\end{figure}

	The isoparameteric transformation is applied to calculate the averages on boundary faces and curved cells. The isoparametric transformation is used to calculate the position and normal direction of the Gaussian integral points on curved mesh interface. The transformation is shown in Figure \ref{fig:quad-curve}, and is written as
	\begin{equation} \label{line-iso}
		\begin{aligned}
			\mathbf{x}(\xi) &= \sum_{k=1}^{3} v_k(\xi) \mathbf{x}_k, \\
			v_1(\xi) &= 1-3\xi+2\xi^2, \\
			v_2(\xi) &= -\xi+2\xi^2, \\
			v_3(\xi) &= 4\xi-4\xi^2.
		\end{aligned}
	\end{equation}
	
	The normal direction at the Gaussian points on the quadratic curve is calculated by
	\begin{equation}
		\begin{aligned}
			n_x &= \sum_{k=1}^{3} v_k^{'}(\xi) y_k, \\
			n_y &= -\sum_{k=1}^{3} v_k^{'}(\xi) x_k.
		\end{aligned} \nonumber
	\end{equation}
	
	The line-averaged value of variable $Q$, such as the $\rho U, \rho V$ at non-slip wall boundary, is calculated by
	\begin{equation} \label{line-aver}
		\begin{aligned}
			\frac{1}{L}\int_l Q \mathbf{d} s
			&= \frac{1}{L} \int_0^1 Q(x(\xi), y(\xi)) \sqrt{\dot{x}^2 + \dot{y}^2} \mathbf{d} \xi, \\
		\end{aligned}
	\end{equation}
	where $L$ is the arc-length of the quadratic curve. Similarly, the line-averaged normal gradients of variable $Q$ is calculated by
	\begin{equation} \label{line-grad-aver}
		\begin{aligned}
			\frac{1}{L}\int_l \frac{\partial Q}{\partial \mathbf{n}} \mathbf{d} s
			= \frac{1}{L}\int_0^1 (\frac{\partial Q}{\partial x}n_x + \frac{\partial Q}{\partial y}n_y)  \sqrt{\dot{x}^2 + \dot{y}^2} \mathbf{d} \xi.
		\end{aligned}
	\end{equation}
	
\begin{figure}[!htbp]
		\centering
		\includegraphics[width = 0.75\columnwidth]{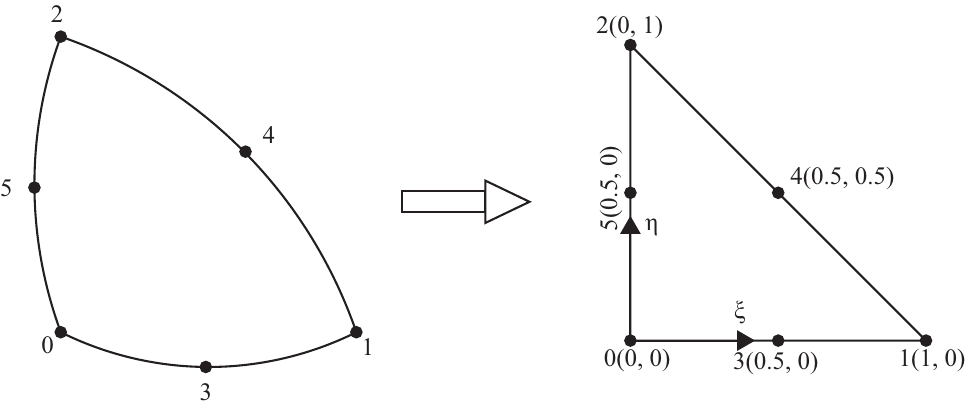}
		\caption{Isoparametric transformation for quadratic triangular cell}
		\label{fig:map-tri}
	\end{figure}

	The integral in Eq. \ref{line-aver} and \ref{line-grad-aver} are calculated by its exact solution, as shown in Appendix A.
	The isoparameteric transformation for quadratic triangular cell is shown in Figure \ref{fig:map-tri}, and the transformation is written as
	\begin{equation} \label{tri-iso}
		\begin{aligned}
			\mathbf{x}(\xi, \eta) &= \sum_{k=1}^{6} v_k(\xi, \eta) \mathbf{x}_k, \\
			v_1(\xi, \eta) &= (\xi + \eta - 1)(2\xi + 2\eta - 1), \\
			v_2(\xi, \eta) &= \xi (2\xi - 1), \\
			v_3(\xi, \eta) &= \eta (2\eta - 1), \\
			v_4(\xi, \eta) &= -4\xi (\xi + \eta - 1), \\
			v_5(\xi, \eta) &= 4\xi \eta, \\
			v_6(\xi, \eta) &= -4\eta (\xi + \eta - 1).
		\end{aligned}
	\end{equation}
	
	One curved quadrilateral cell is divided into two curved triangles to calculate the cell averages. The cell-averages of the basis functions over the curved triangular cell in Eq. \ref{p_2} is calculated by
	\begin{equation} \label{cell_aver}
		\iint_{\Omega} \phi_k(x, y) {\rm d}x{\rm d}y = \iint_{\tilde{\Omega}} \phi_k(x(\xi, \eta), y(\xi, \eta)) J(\xi, \eta) {\rm d}\xi{\rm d}\eta,
	\end{equation}
	where $\Omega$ indicates the grid cell in physical coordinate system, $\tilde{\Omega}$ indicates the corresponding integral domain in reference coordinate system, $J$ is the Jocabian matrix of the isoparametric transformation, $\xi, \eta$ are the axes in reference coordinate system. The integral in Eq. \ref{cell_aver} is also calculated by its exact solution, as shown in Appendix B.
	
	\subsection{Error analysis} \label{err-ana-sec}
\begin{figure}[!htbp]
		\centering
		\includegraphics[width = 0.5\columnwidth]{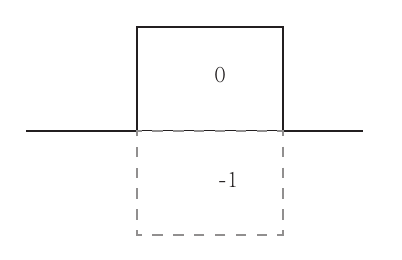}
		\caption{Error analysis of ghost cell method}
		\label{fig:error_ana}
	\end{figure}

	In this section, the error introduced by the setting of ghost states and the proposed one-sided reconstruction is analyzed in a simplified model as shown in Figure \ref{fig:error_ana}.

	To constrain the flow variables on wall boundary, for example the velocity on non-slip walls, the flow variables in the ghost cell is set as $\overline{Q}_{-1} = -\overline{Q}_0$, with the condition
	\begin{equation}
		\overline{Q}_0 - a_2 \Delta y+ a_5 \Delta y^2 = \overline{Q}_{-1}, \nonumber
	\end{equation}
	which is equivalent to,
	\begin{equation} \label{ghost-u}
		\overline{Q}_0 - \frac{1}{2} a_2 \Delta y+ \frac{1}{2} a_5 \Delta y^2 = 0.
	\end{equation}
	
	The flow variable on boundary is directly restricted in the one-sided reconstruction, which provides a condition
	\begin{equation}  \label{bn-u}
		\overline{Q}_0 - \frac{1}{2} a_2 \Delta y+ \frac{1}{6} a_5 \Delta y^2 = 0.
	\end{equation}
	
	Considering that the conditions Eq. \ref{ghost-u} and Eq. \ref{bn-u} are constrained in constrained least-square method, the setting of ghost state will introduce an $o(\Delta y^2)$ error at the wall boundary. To eliminate this error, high-order extrapolation methods for the construction of ghost state, for example the inverse Lax-Wendroff method \cite{invLW}, are needed. But the third-order scheme with ghost cells still has better boundary layer resolution than the second-order treatment, because there are more degrees of freedom in each cell. When the boundary layer is ill-resolved by the grid, the reconstructed velocity on boundary by ghost cell method is not small enough, which lead to slip flow near the wall boundary.
	
	To constrain the normal gradient of the flow variable on boundary, for example the pressure on non-slip walls, the flow variables in the ghost cell is set as $\overline{Q}_{-1} = \overline{Q}_0$, and
	\begin{equation}
		\overline{Q}_0 - a_2 \Delta y+ a_5 \Delta y^2 = \overline{Q}_{-1}, \nonumber
	\end{equation}
	which is equivalent to,
	\begin{equation} \label{ghost-u-grad}
		- a_2 \Delta y+ a_5 \Delta y^2 = 0.
	\end{equation}
	In the one-sided reconstruction, the normal gradients of flow variables on boundary are restricted directly,
	\begin{equation}  \label{bn-u-grad}
		a_2 \Delta y - a_5 \Delta y^2 = 0,
	\end{equation}
	which is identical with Eq. \ref{ghost-u-grad}, and the ghost state $\overline{Q}_{-1} = \overline{Q}_0$ can give a zero-normal derivative.
	
	Above all, the third-order scheme with ghost cells can provide accurate zero-normal gradient and introduce an $o(\Delta y^2)$ error on the variables at the wall boundary.
	
	\section{Numerical Examples}
	In this section, laminar boundary layer, subsonic and hypersonic viscous flow around a cylinder are studied by three types of boundary treatment, namely Green-Gauss scheme on linear grids (Method \Rmnum{1}), third-order compact scheme with ghost cells on linear grids (Method \Rmnum{2}) and the proposed one-sided third-order compact scheme (Method \Rmnum{3}) on curved grids.
	It is noted that the interior cells are reconstructed by the same third-order compact GKS, the only difference is the reconstruction on boundary cells.
	\subsection{Accuracy test for linear scheme}
	For Euler equations, sine wave propagation problem is widely used to test the accuracy of numerical method. The exact solution of this problem is
	\begin{equation}
		\begin{gathered}
			\rho = 1.0 + 0.2 \cdot sin(2\pi(x-t)) \cdot sin(2\pi(y-t)), \\
			U = 1.0, V = 1.0,p = 1.0.
		\end{gathered} \nonumber
	\end{equation}
	The computational domain is $x \in [0,1], y \in [0,1]$. Third-order compact gas-kinetic scheme and S2O4 temporal scheme are adopted. Periodic boundary condition is imposed, and CFL number is taken as 0.5, the output time is $t = 1$. The accuracy test results are shown in Table \ref{acc-test},
	\begin{table}[!htbp]
		\small
		\begin{center}
			\def\temptablewidth{1.0\textwidth}
			{\rule{\temptablewidth}{1pt}}
			\begin{tabular*}{\temptablewidth}{@{\extracolsep{\fill}}c c c}
				$h$ & $L_2 \space \rm{error}$  &$\rm{order}$\\
				\hline
				$1/20$ & $2.19\times 10^{-3}$ &-- \\ 	
				$1/40$ & $2.79\times 10^{-4}$ &2.97\\ 	
				$1/80$ & $3.50\times 10^{-5}$ &2.99\\
				$1/160$ & $4.38\times 10^{-6}$ &3.00\\
			\end{tabular*}
			{\rule{\temptablewidth}{1pt}}
		\end{center}
		\vspace{-4mm} \caption{ Accuracy test for linearly third-order scheme}
		\label{acc-test}
	\end{table}
	\subsection{Laminar boundary layer}
\begin{figure}[!htbp]
		\centering
		\subfigure[Geometry]{\includegraphics[width = 0.45\columnwidth]{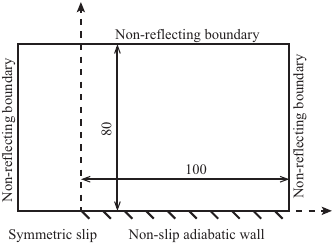}} \hspace{5pt}
		\subfigure[Mesh]{\includegraphics[width = 0.45\columnwidth, trim = 0 0 0 100, clip]{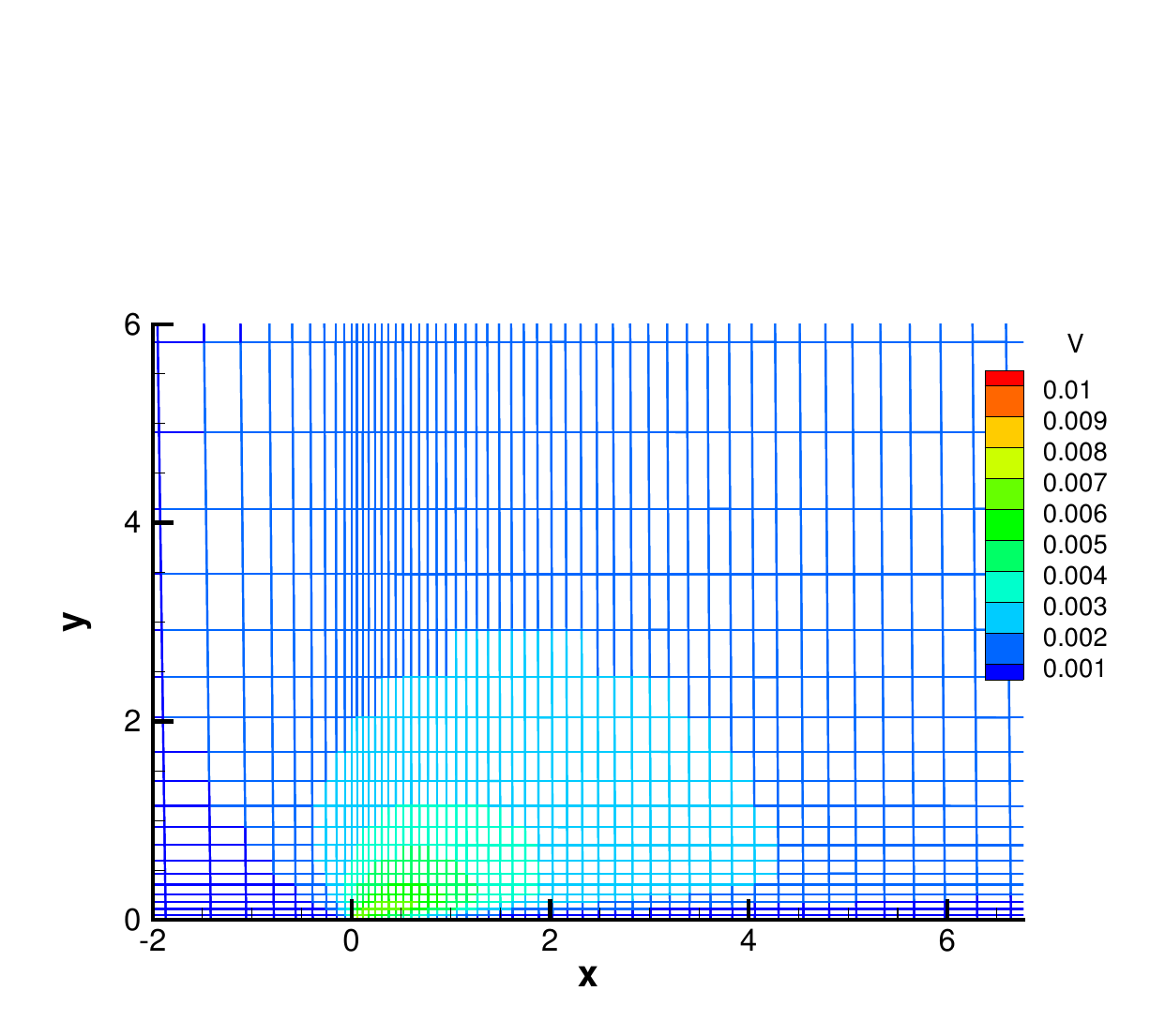}} \\
		\caption{Computational domain of flat plate}
		\label{bnLayer}
	\end{figure}

	A laminar boundary layer over a flat plate is simulated to test the performance of these three boundary treatments on linear mesh. The incoming Mach number $\mathrm{Ma} = 0.15$, Reynolds number $\mathrm{Re} = U_{\infty} L / v = 10^5$ and the characteristic length $L = 100$. The computational domain is shown in Figure \ref{bnLayer}, where the flat plate is placed at $x>0, y=0$. The computation domain, $[-30, 100] \times [0, 80]$, is discretized by $120 \times 35$ quadrilateral cells with the first layer height $h=0.05$ (Mesh \Rmnum{1}), $h=0.1$ (Mesh \Rmnum{2}), $h=0.2$ (Mesh \Rmnum{3}).
	
	The adiabatic non-slip boundary condition is imposed on the plate and symmetric slip boundary condition is set in the front of the plate. The non-reflecting boundary condition based on the Riemann invariants is adopted for the other boundaries, and the incoming flow state is set as $\rho_{\infty} = 1, p_{\infty} = 1/\gamma$.
	Since the flow is rather smooth, the smooth GKS solver Eq. \ref{gks-smooth} and linear compact third-order scheme are adopted to reduce the numerical dissipation, and CFL number is set as $\mathrm{CFL} = 0.3$.\\
	
\begin{figure}[!htbp]
		\centering
		\subfigure{\includegraphics[width = 0.45\columnwidth, trim = 0 10 0 10, clip]{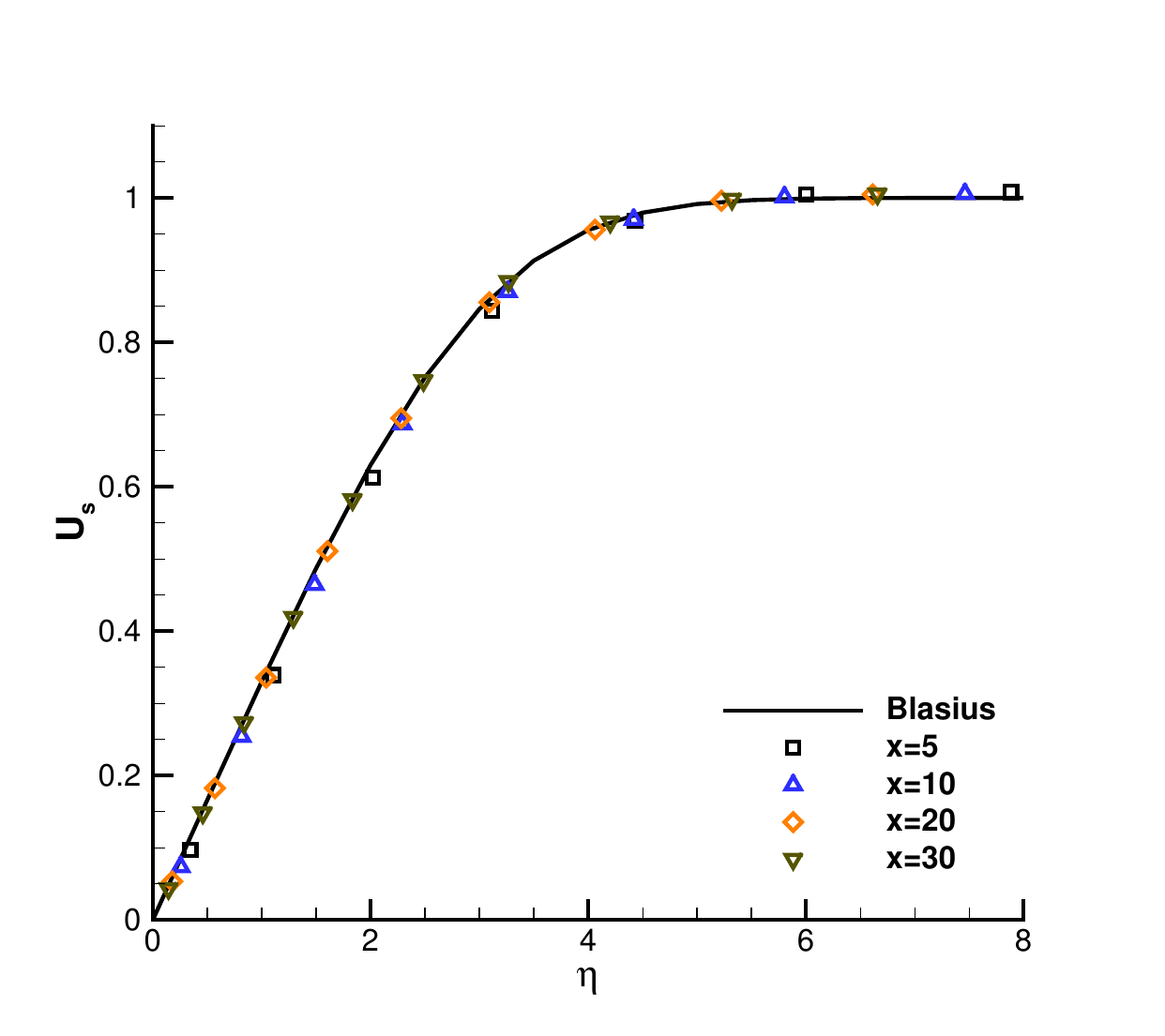}} \hspace{5pt}
		\subfigure{\includegraphics[width = 0.45\columnwidth, trim = 0 10 0 10, clip]{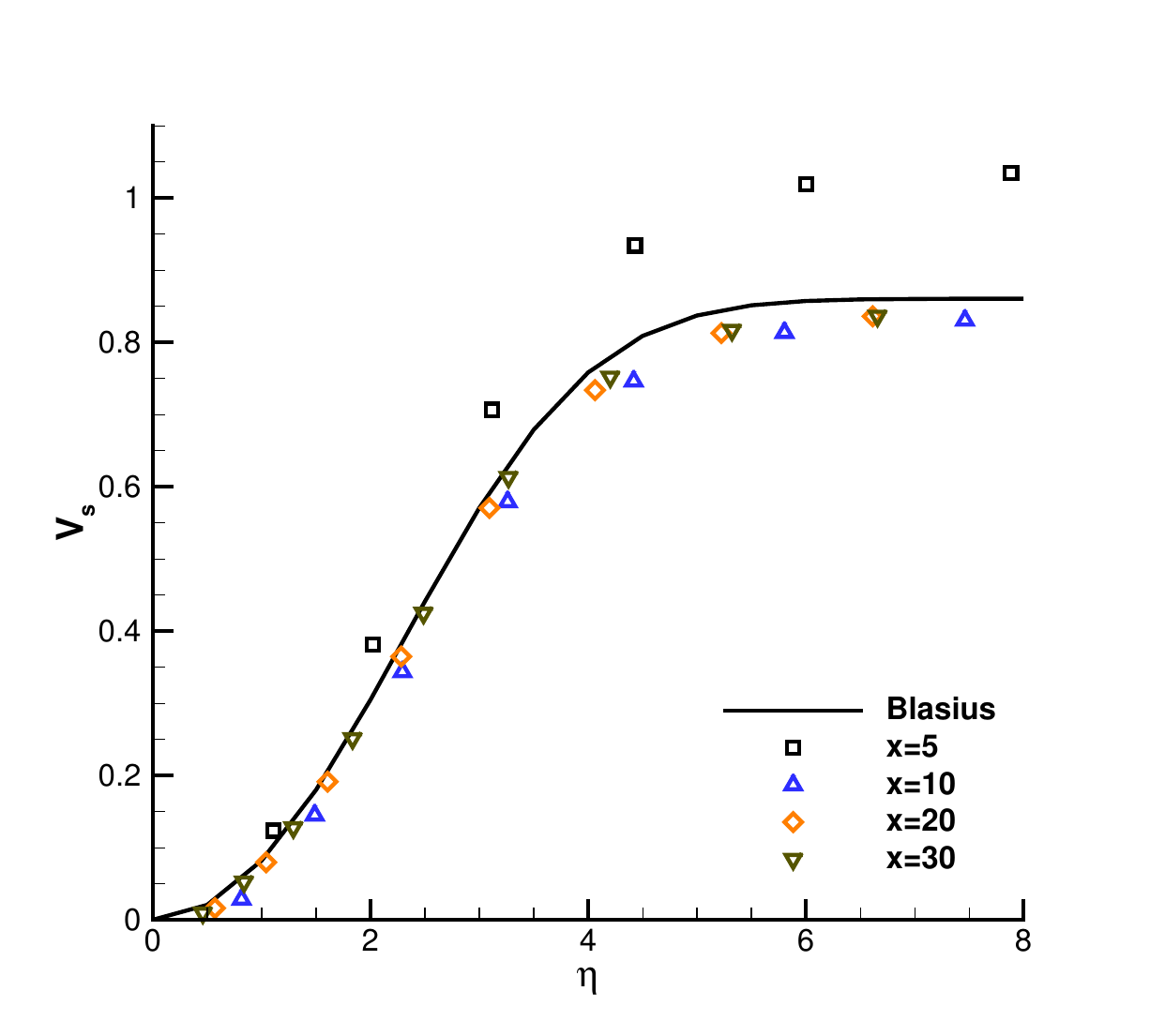}} \\
		\subfigure{\includegraphics[width = 0.45\columnwidth, trim = 0 10 0 10, clip]{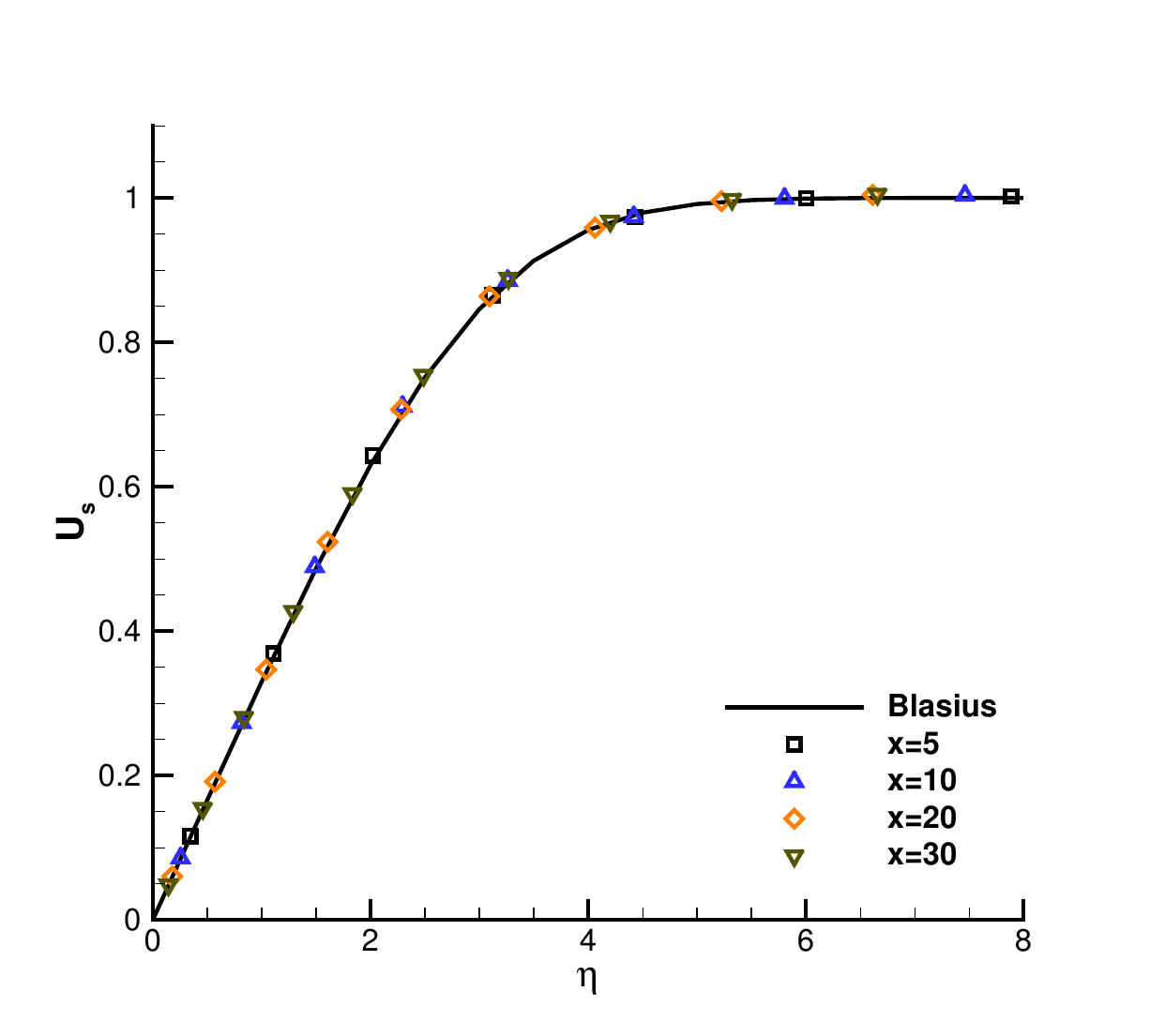}} \hspace{5pt}
		\subfigure{\includegraphics[width = 0.45\columnwidth, trim = 0 10 0 10, clip]{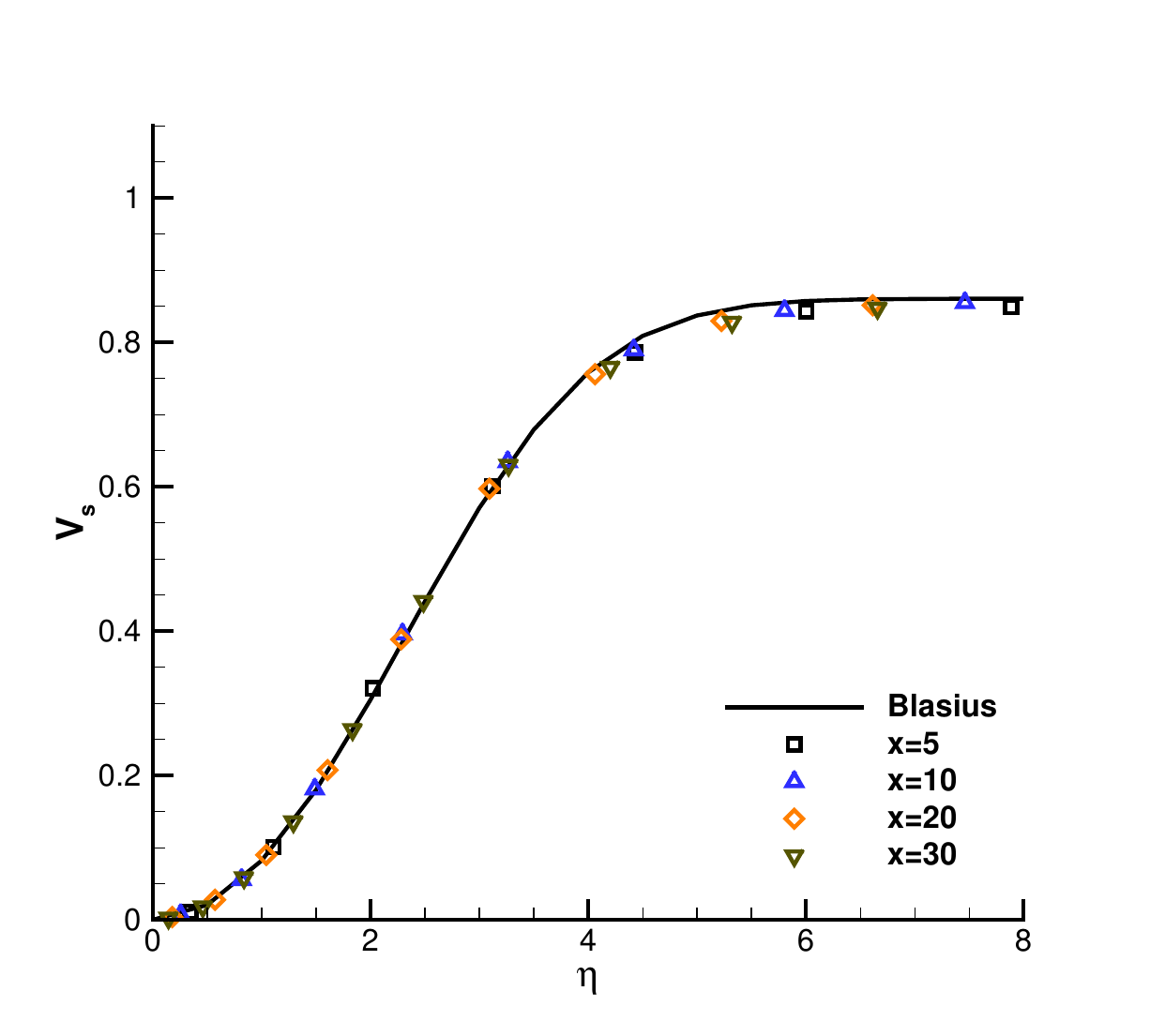}} \\
		\subfigure{\includegraphics[width = 0.45\columnwidth, trim = 0 10 0 10, clip]{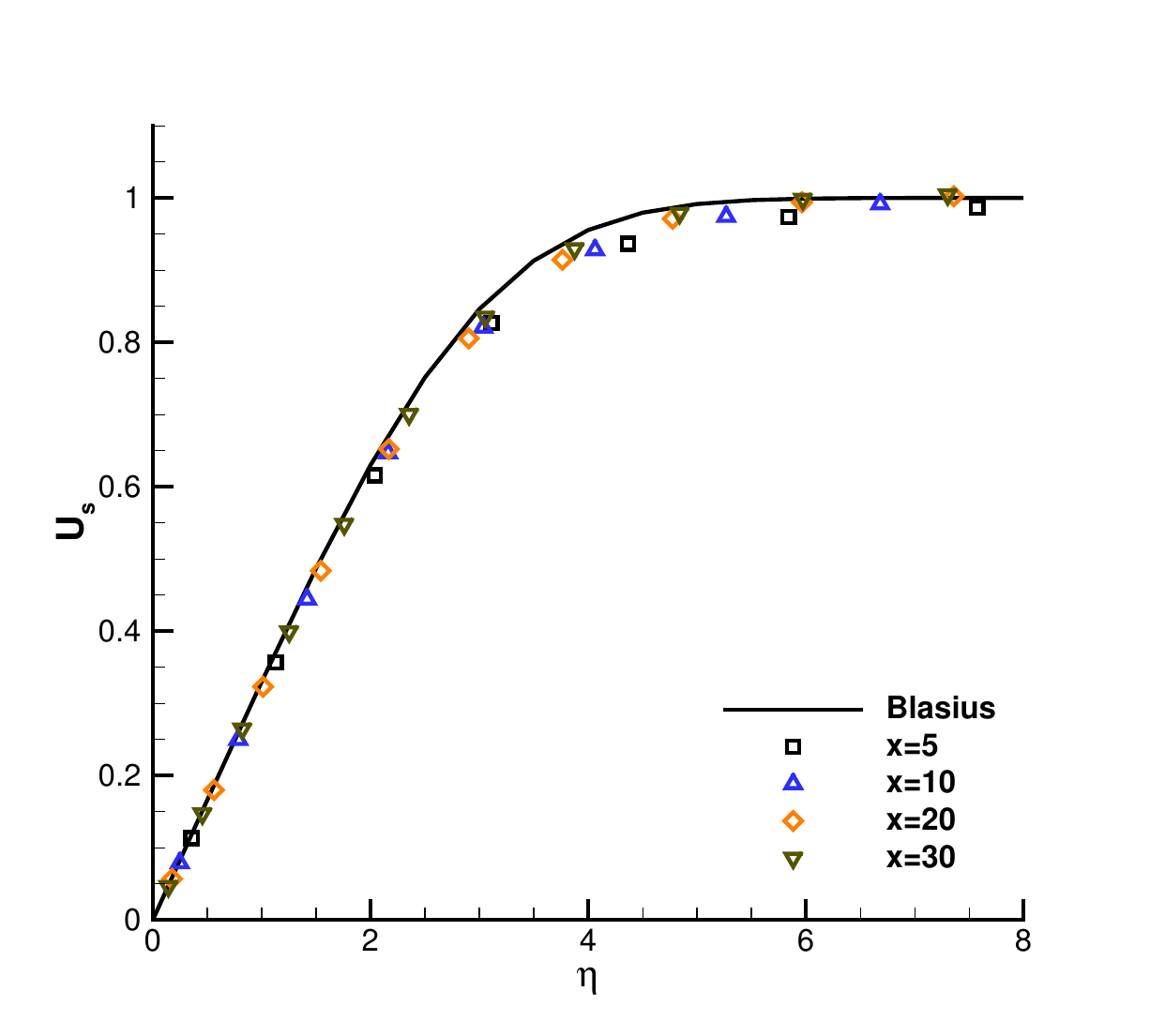}} \hspace{5pt}
		\subfigure{\includegraphics[width = 0.45\columnwidth, trim = 0 10 0 10, clip]{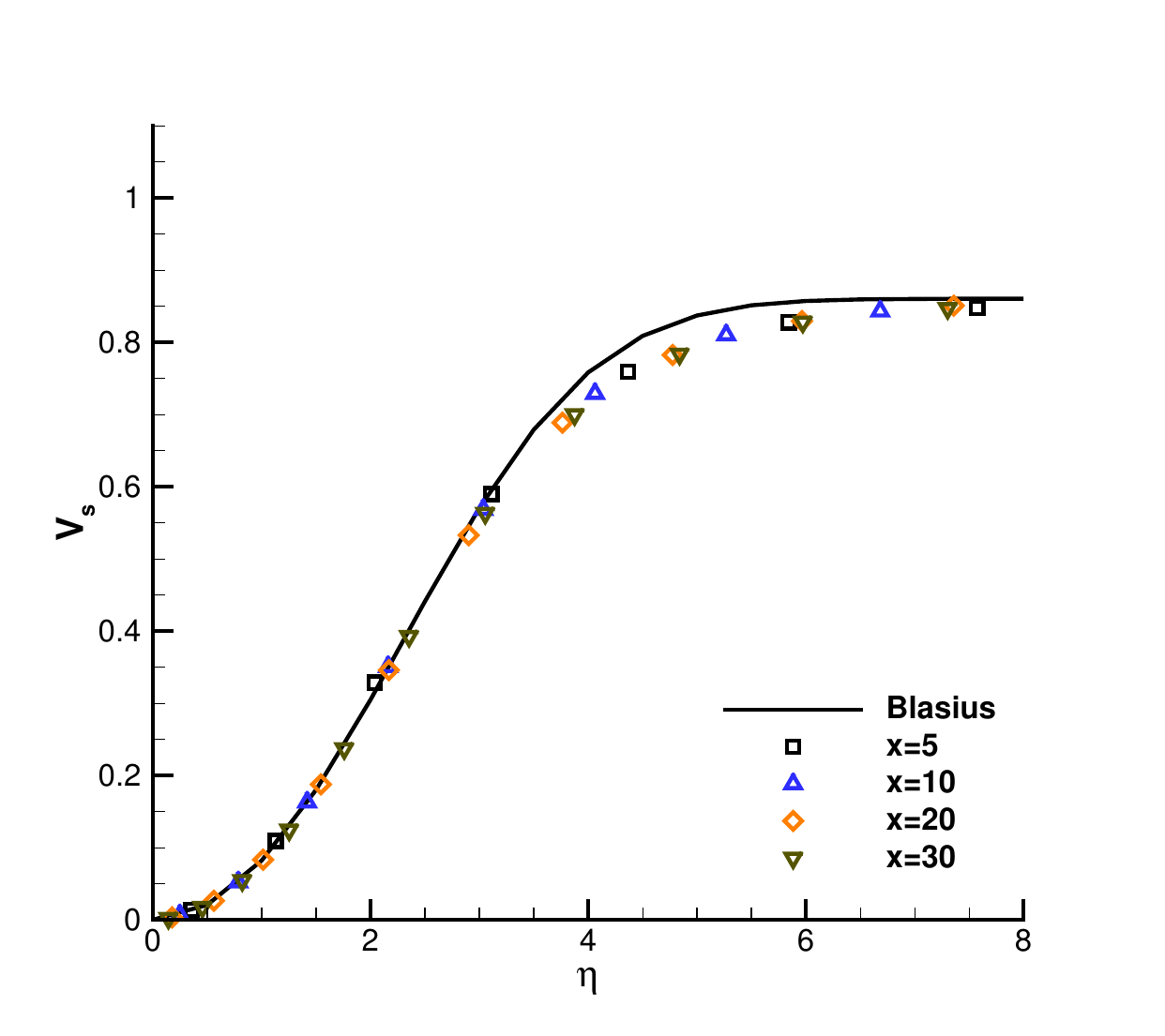}} \\
		\caption{Velocity profile on Mesh \Rmnum{1} with different boundary treatment methods(left column: $U_s$ profile, right column: $V_s$ profile, top row: Method \Rmnum{1}, middle row: Method \Rmnum{2}, bottom row: Method \Rmnum{3})}
		\label{h005}
	\end{figure}

\begin{figure}[!htbp]
		\centering
		\subfigure{\includegraphics[width = 0.45\columnwidth, trim = 0 10 0 10, clip]{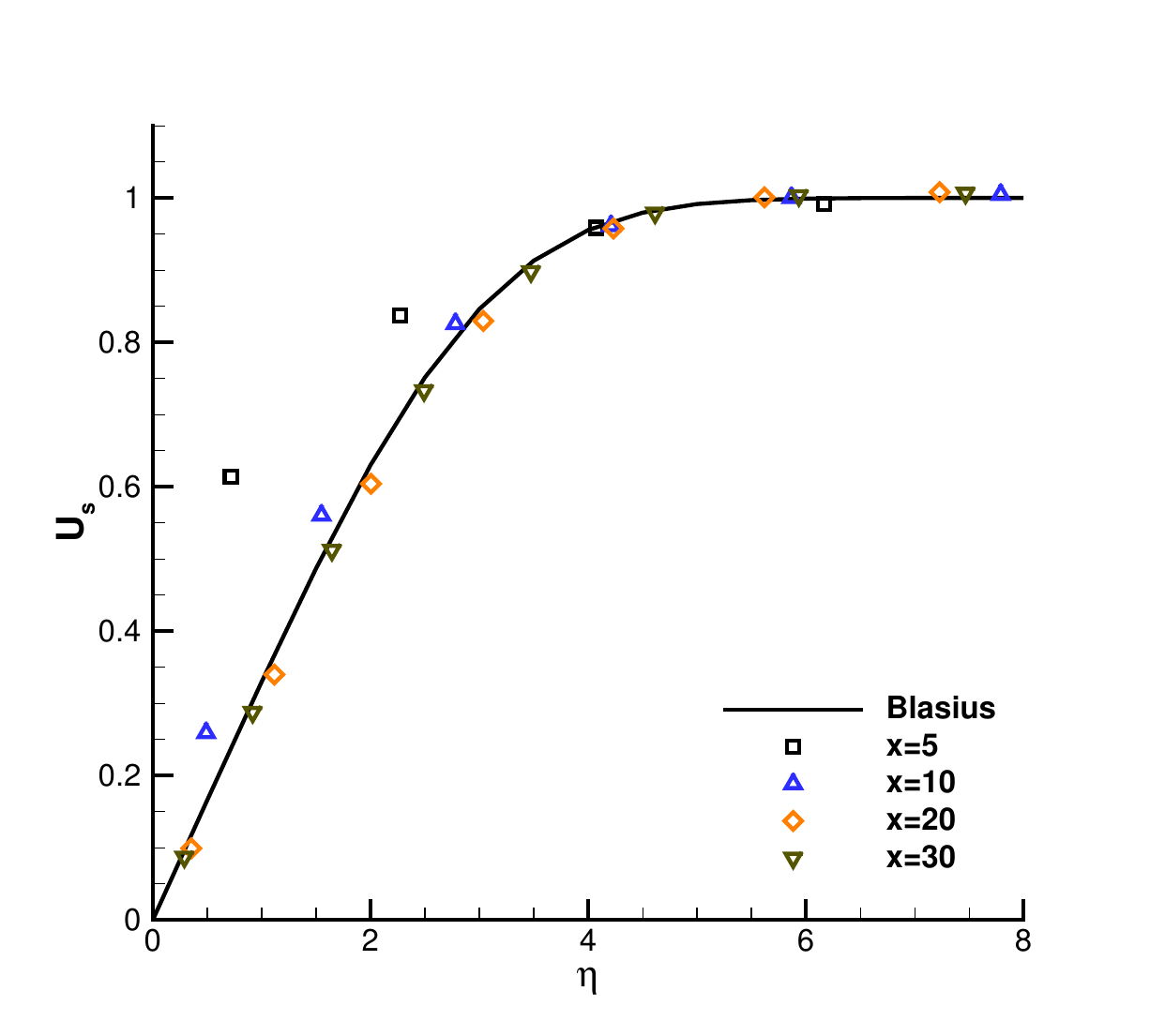}} \hspace{5pt}
		\subfigure{\includegraphics[width = 0.45\columnwidth, trim = 0 10 0 10, clip]{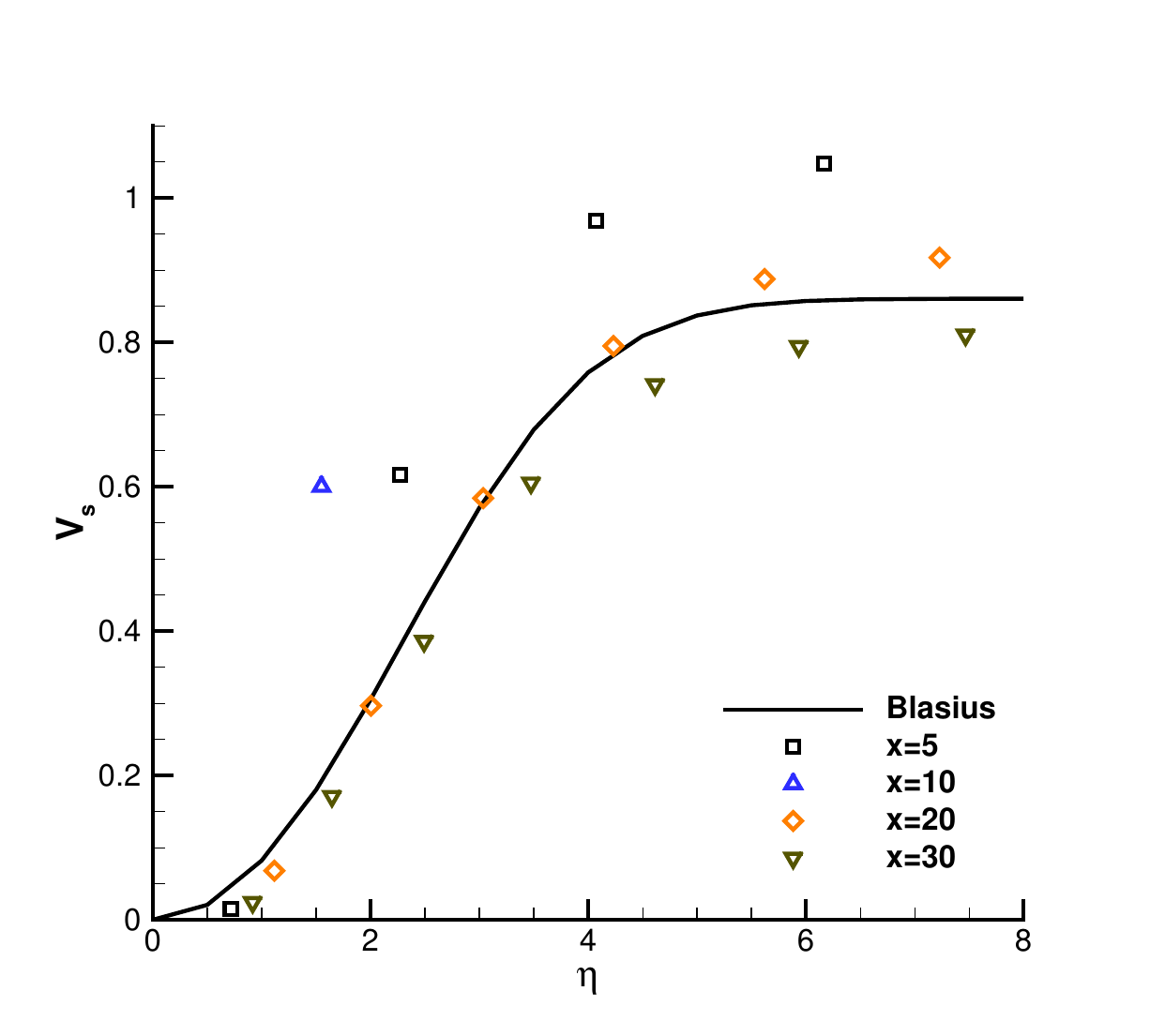}} \\
		\subfigure{\includegraphics[width = 0.45\columnwidth, trim = 0 10 0 10, clip]{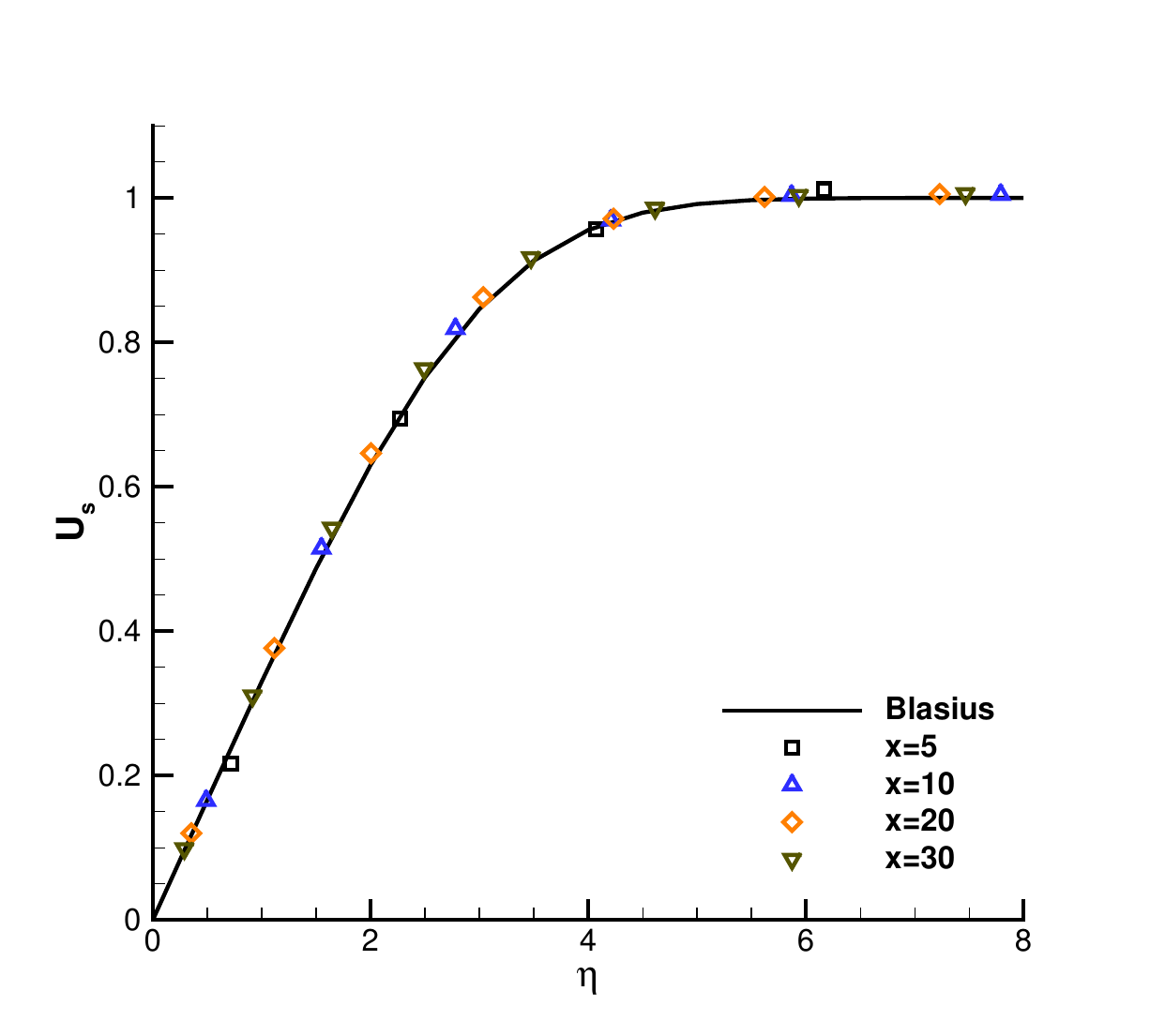}} \hspace{5pt}
		\subfigure{\includegraphics[width = 0.45\columnwidth, trim = 0 10 0 10, clip]{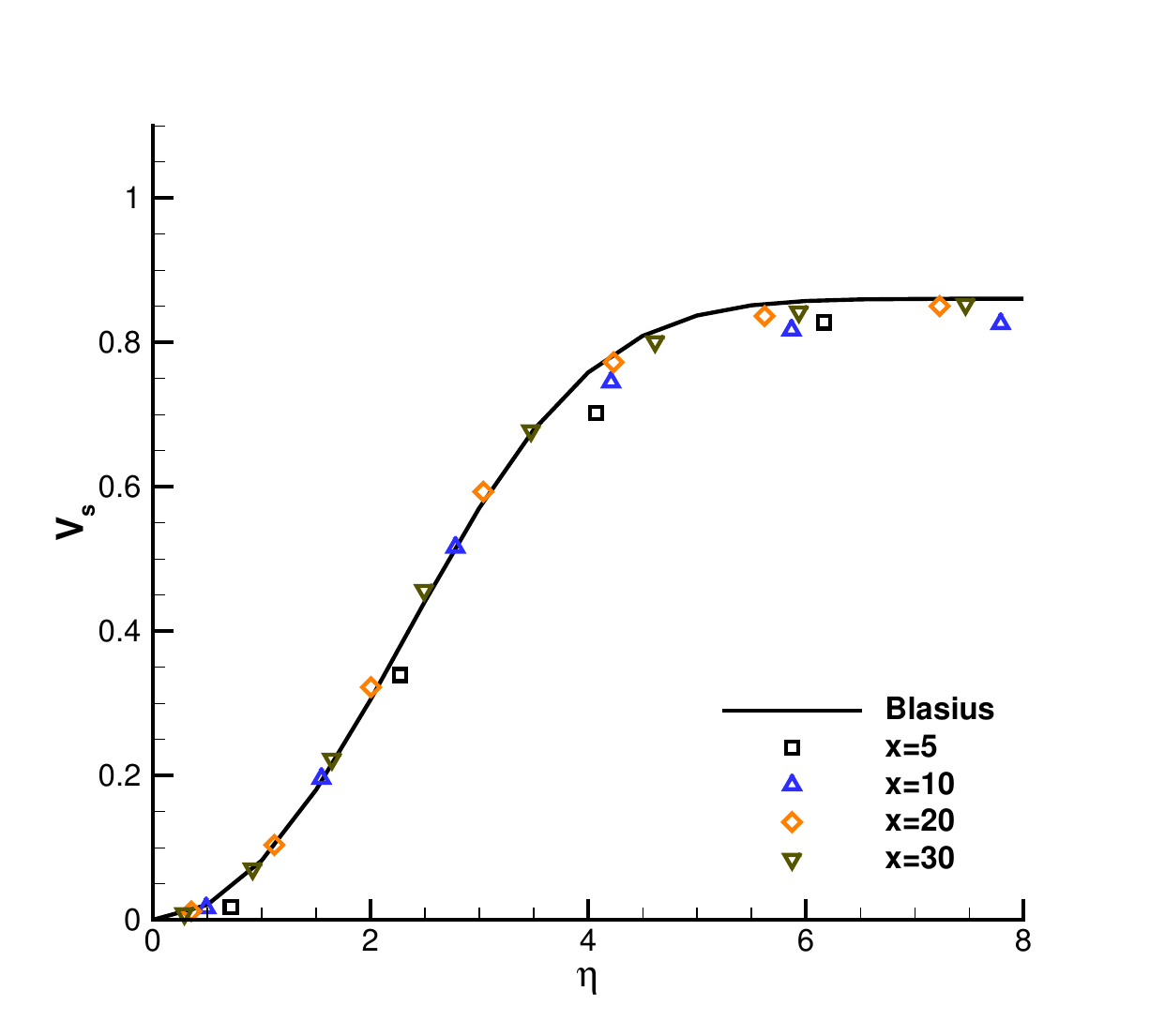}} \\
		\subfigure{\includegraphics[width = 0.45\columnwidth, trim = 0 10 0 10, clip]{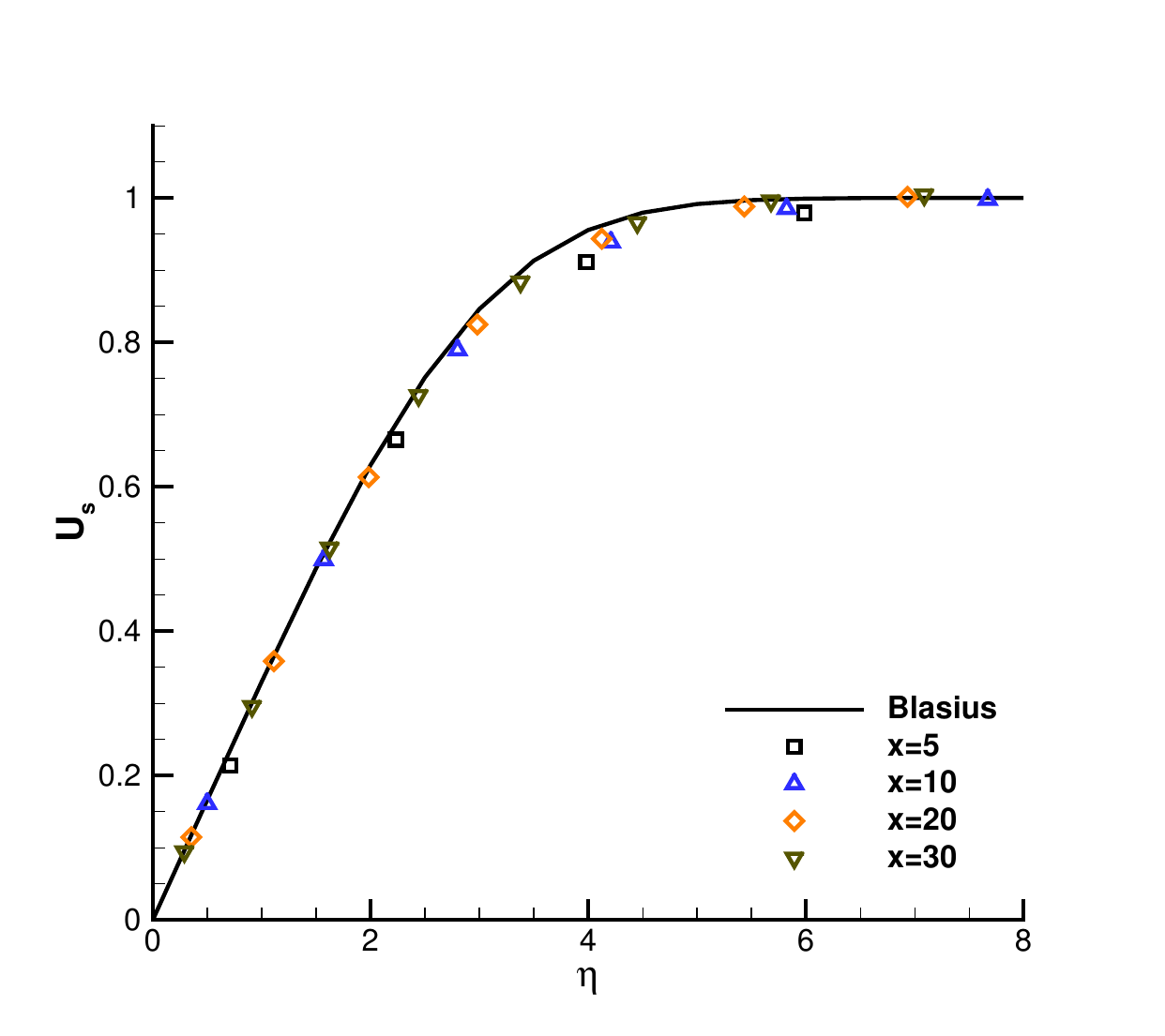}} \hspace{5pt}
		\subfigure{\includegraphics[width = 0.45\columnwidth, trim = 0 10 0 10, clip]{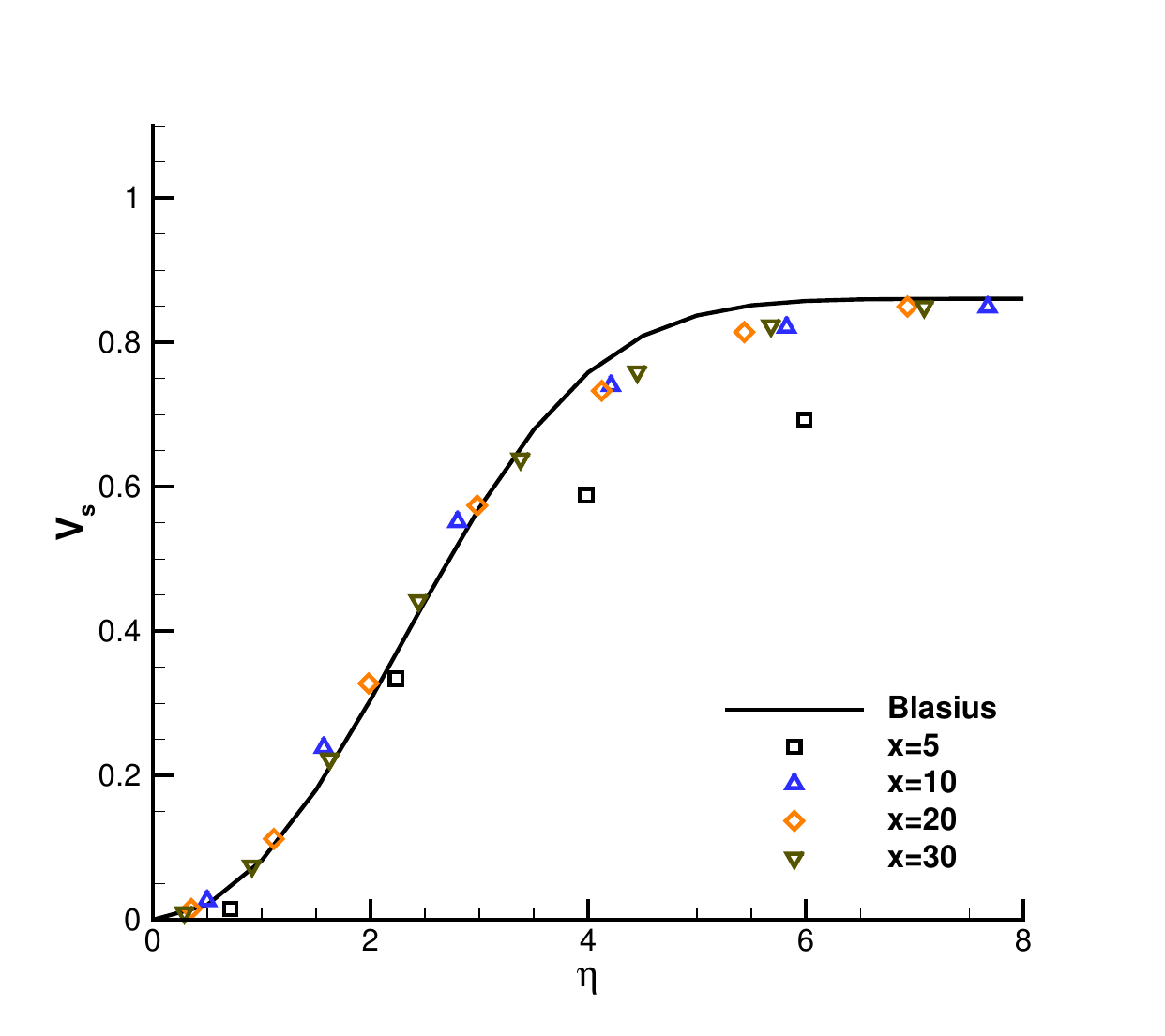}} \\
		\caption{Velocity profile on Mesh \Rmnum{2} with different boundary treatment methods(left column: $U_s$ profile, right column: $V_s$ profile, top row: Method \Rmnum{1}, middle row: Method \Rmnum{2}, bottom row: Method \Rmnum{3})}
		\label{h010}
	\end{figure}

\begin{figure}[!htbp]
		\centering
		\subfigure{\includegraphics[width = 0.45\columnwidth, trim = 0 10 0 10, clip]{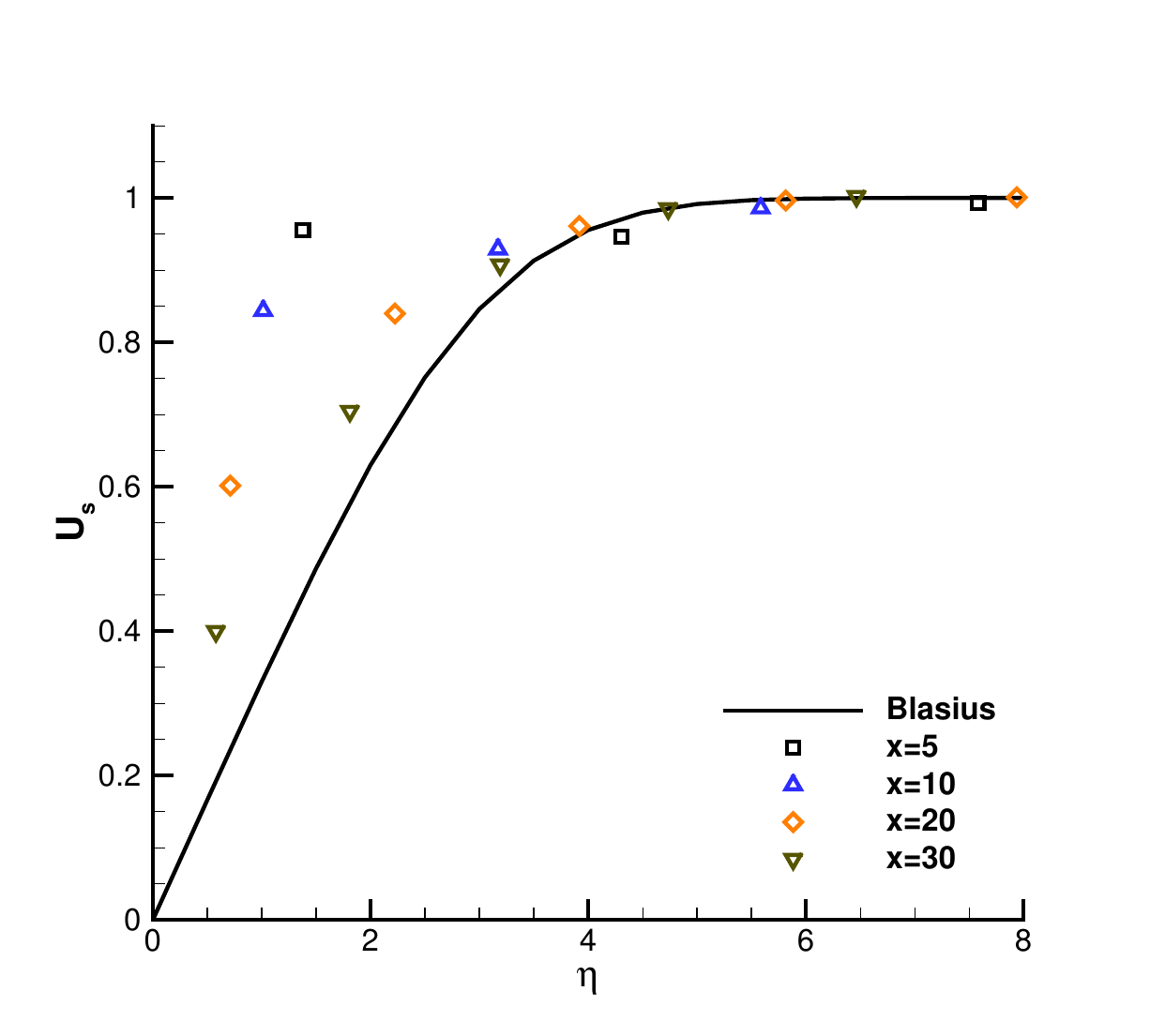}} \hspace{5pt}
		\subfigure{\includegraphics[width = 0.45\columnwidth, trim = 0 10 0 10, clip]{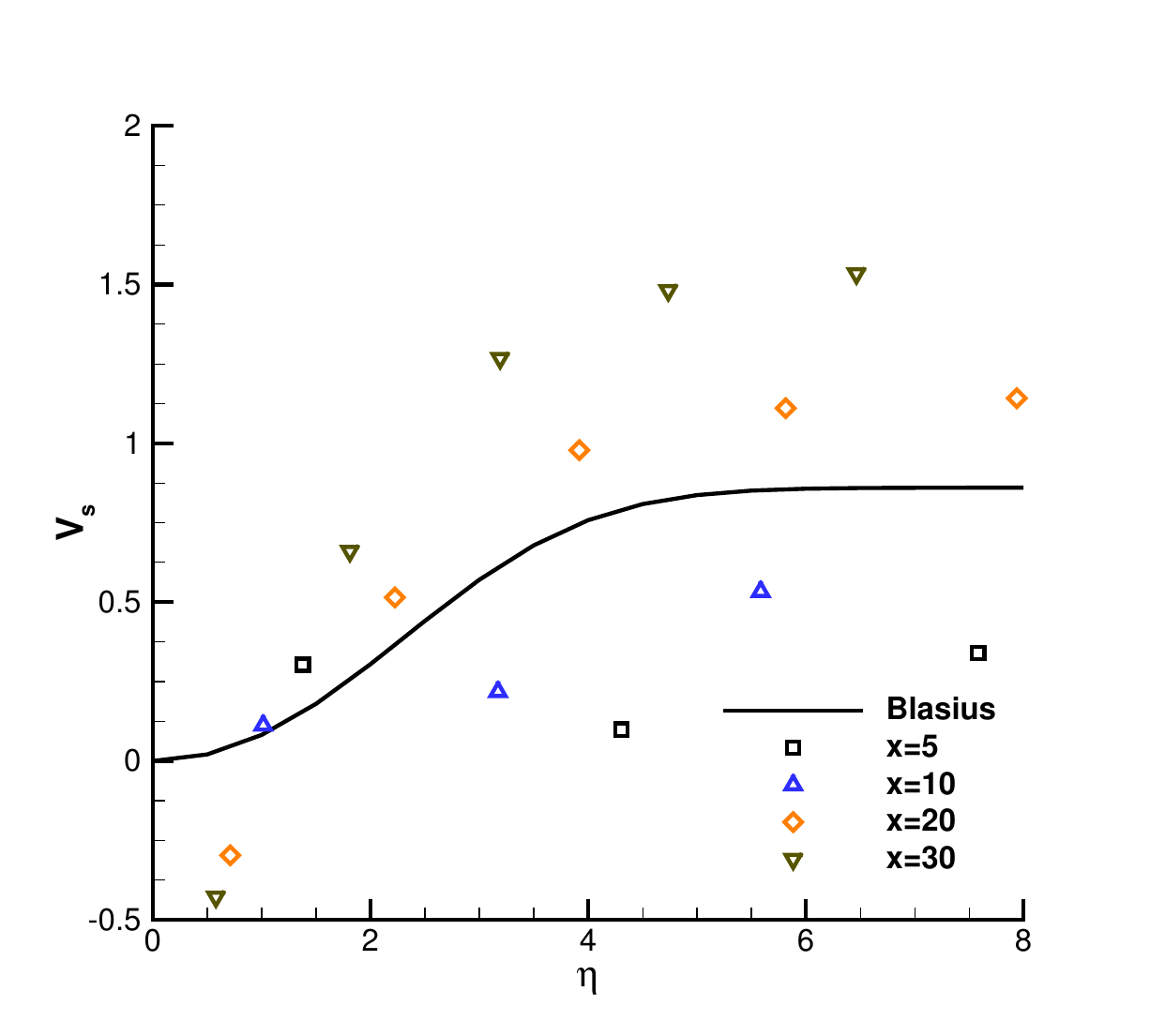}} \\
		\subfigure{\includegraphics[width = 0.45\columnwidth, trim = 0 10 0 10, clip]{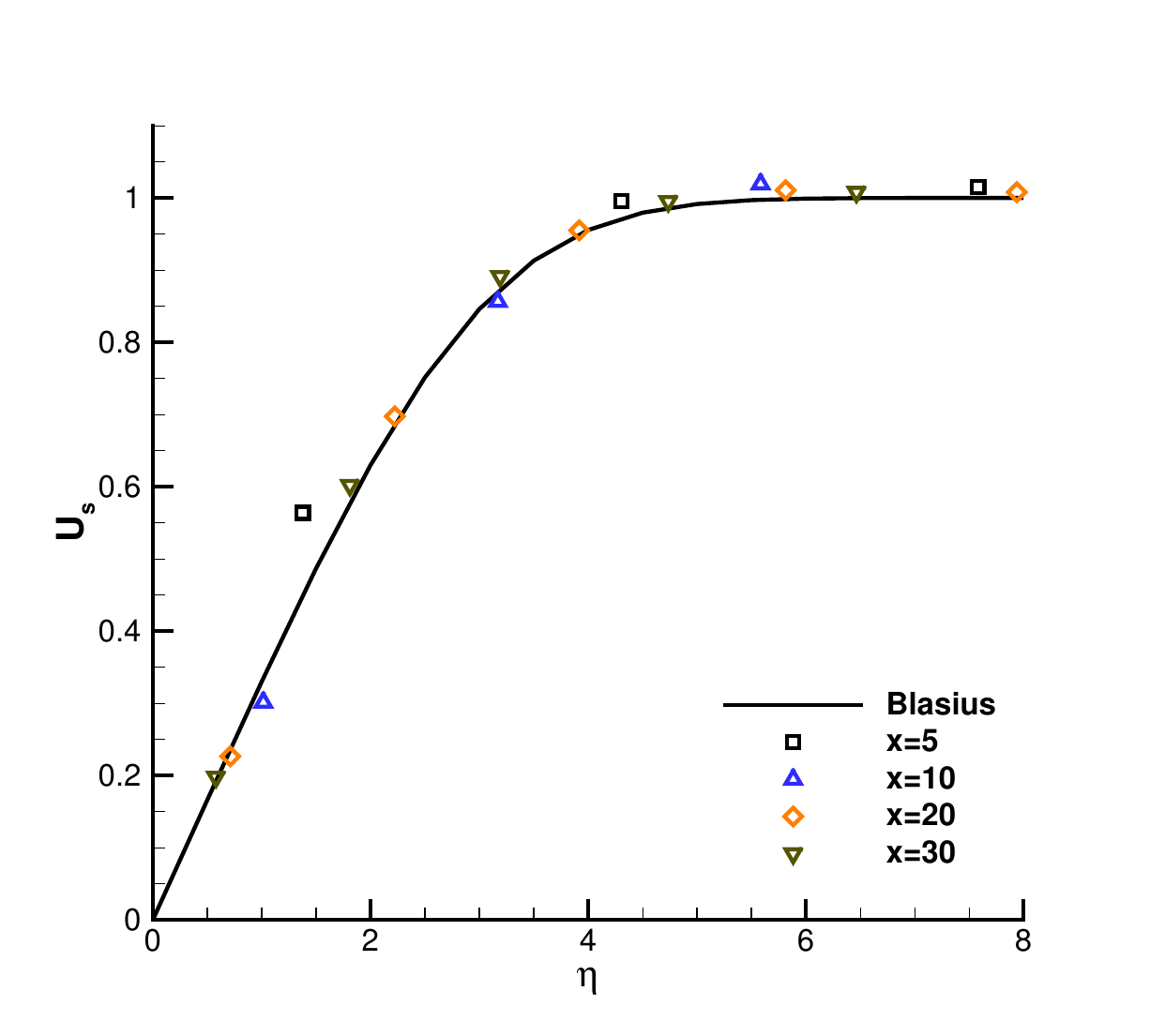}} \hspace{5pt}
		\subfigure{\includegraphics[width = 0.45\columnwidth, trim = 0 10 0 10, clip]{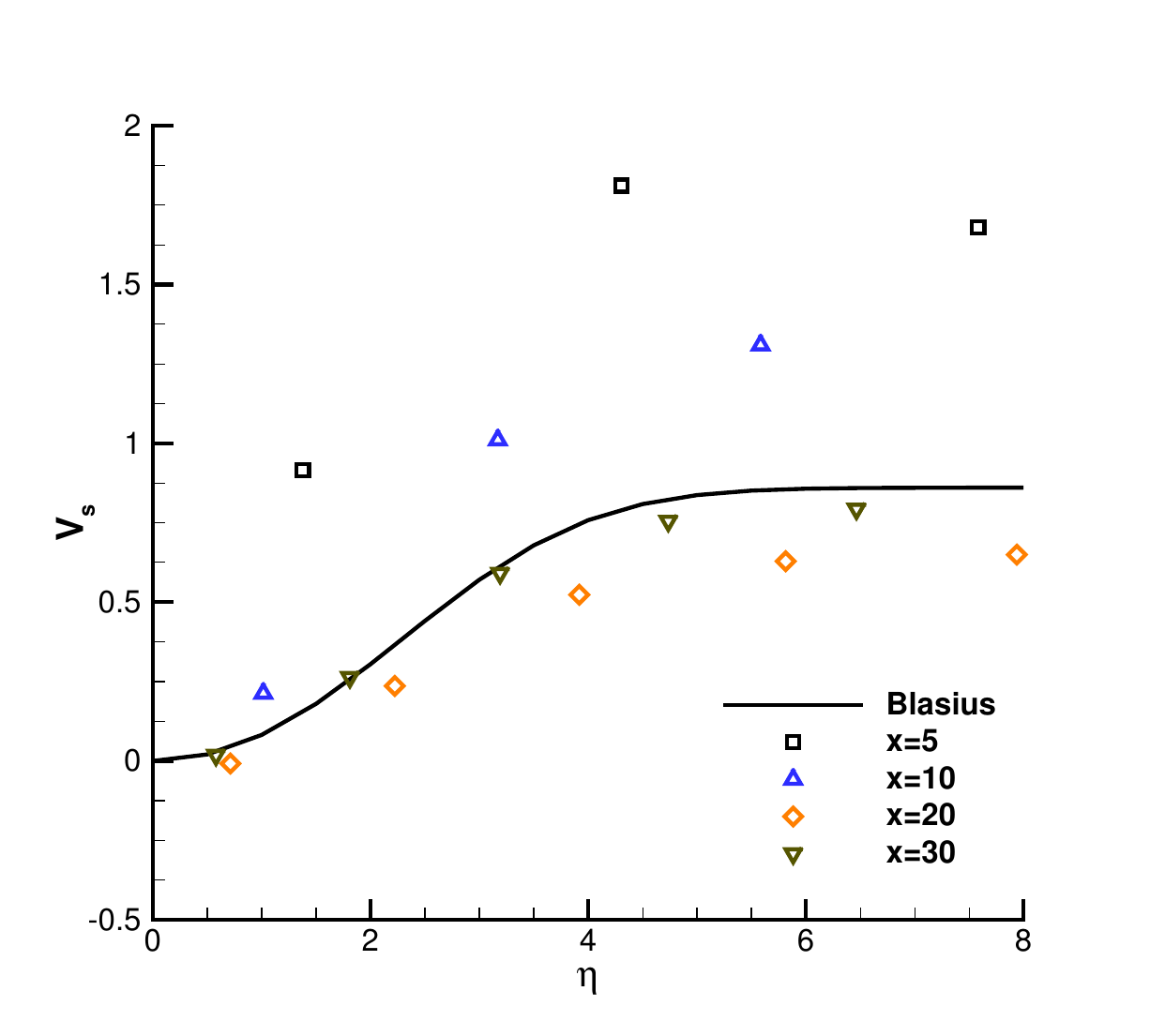}} \\
		\subfigure{\includegraphics[width = 0.45\columnwidth, trim = 0 10 0 10, clip]{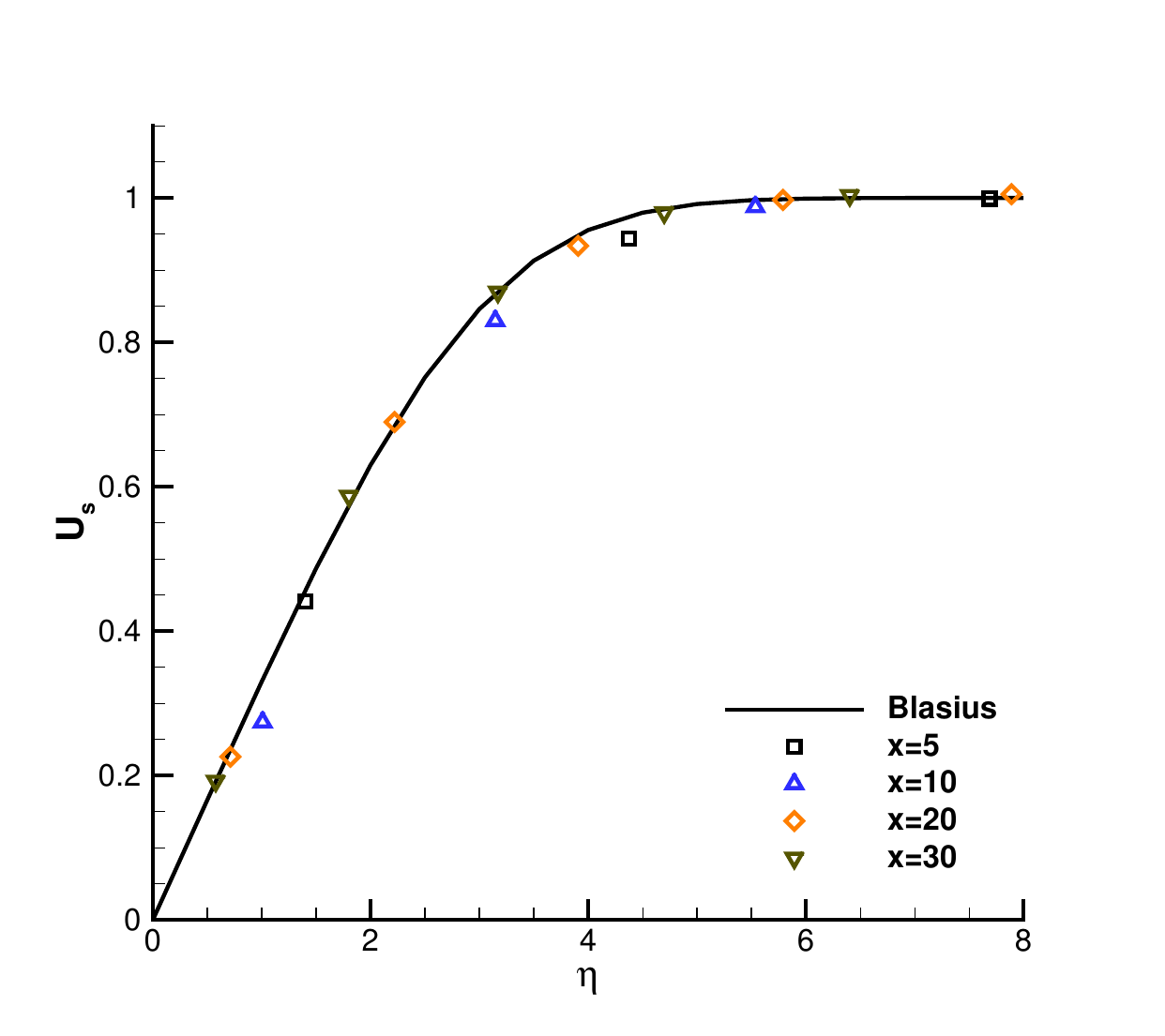}} \hspace{5pt}
		\subfigure{\includegraphics[width = 0.45\columnwidth, trim = 0 10 0 10, clip]{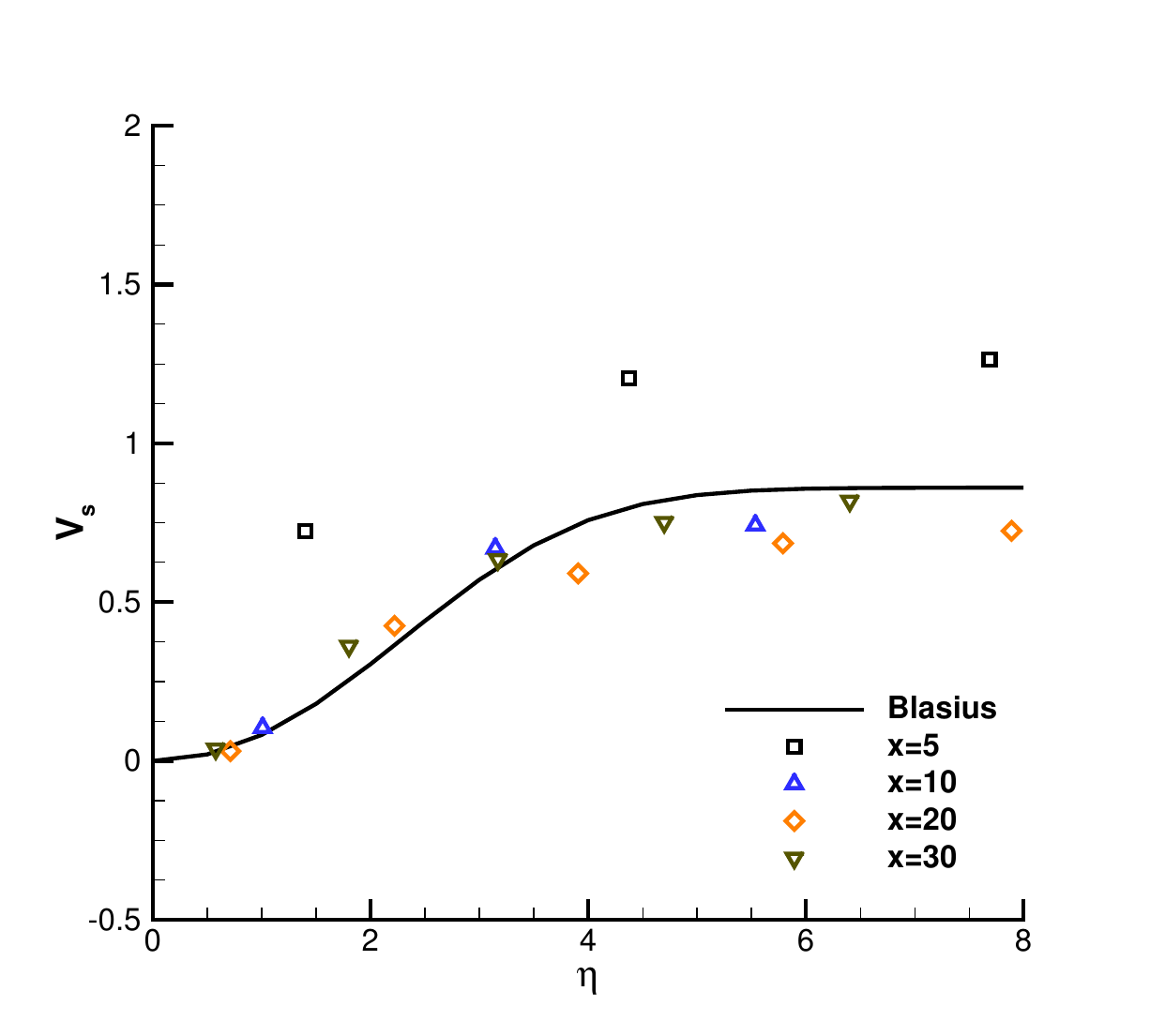}} \\
		\caption{Velocity profile on Mesh \Rmnum{3}  with different boundary treatment methods(left column: $U_s$ profile, right column: $V_s$ profile, top row: Method \Rmnum{1}, middle row: Method \Rmnum{2}, bottom row: Method \Rmnum{3})}
		\label{h020}
	\end{figure}

	The velocity and skin-friction is normalized in the following way. The normalized distance from the flat plate is defined as
	\begin{equation}
		\eta = y \sqrt{\frac{\rho U_{\infty}}{x \mu_{\infty}}}. \nonumber
	\end{equation}
	The velocity is normalized by
	\begin{equation}
		\begin{gathered}
			U_s = \frac{U}{U_{\infty}}, \\
			V_s = V \sqrt{\frac{\rho x}{\mu_{\infty} U_{\infty}}}.
		\end{gathered} \nonumber
	\end{equation}
	The local Reynolds number $\mathrm{Re}_x$ and the skin-friction coefficient is defined as
	\begin{equation}
		\mathrm{Re}_x = \frac{x}{L}\mathrm{Re}, C_f = \frac{2\tau_{wall}}{\rho_{\infty} U_{\infty}^2}, \nonumber
	\end{equation}
	where $\tau_{wall}$ is the skin shear stress.
	The profiles of normalized velocity at $x = 5, 10, 20, 30$ are shown in Figure \ref{h005}, \ref{h010} and \ref{h020}. It is shown that the two third-order methods, Method \Rmnum{2} and Method \Rmnum{3} can achieve much better performance than the second-order Method \Rmnum{1}. Even the results of Method \Rmnum{2} and Method \Rmnum{3} on Mesh \Rmnum{2} are better than those of Method \Rmnum{1} on Mesh \Rmnum{1}. And the proposed Method \Rmnum{3} performs better than Method \Rmnum{2} on Mesh \Rmnum{3}, which is consistent with the analysis in Section \ref{err-ana-sec}.
	
\begin{figure}[!htbp]
		\centering
		\includegraphics[width = 0.5\columnwidth]{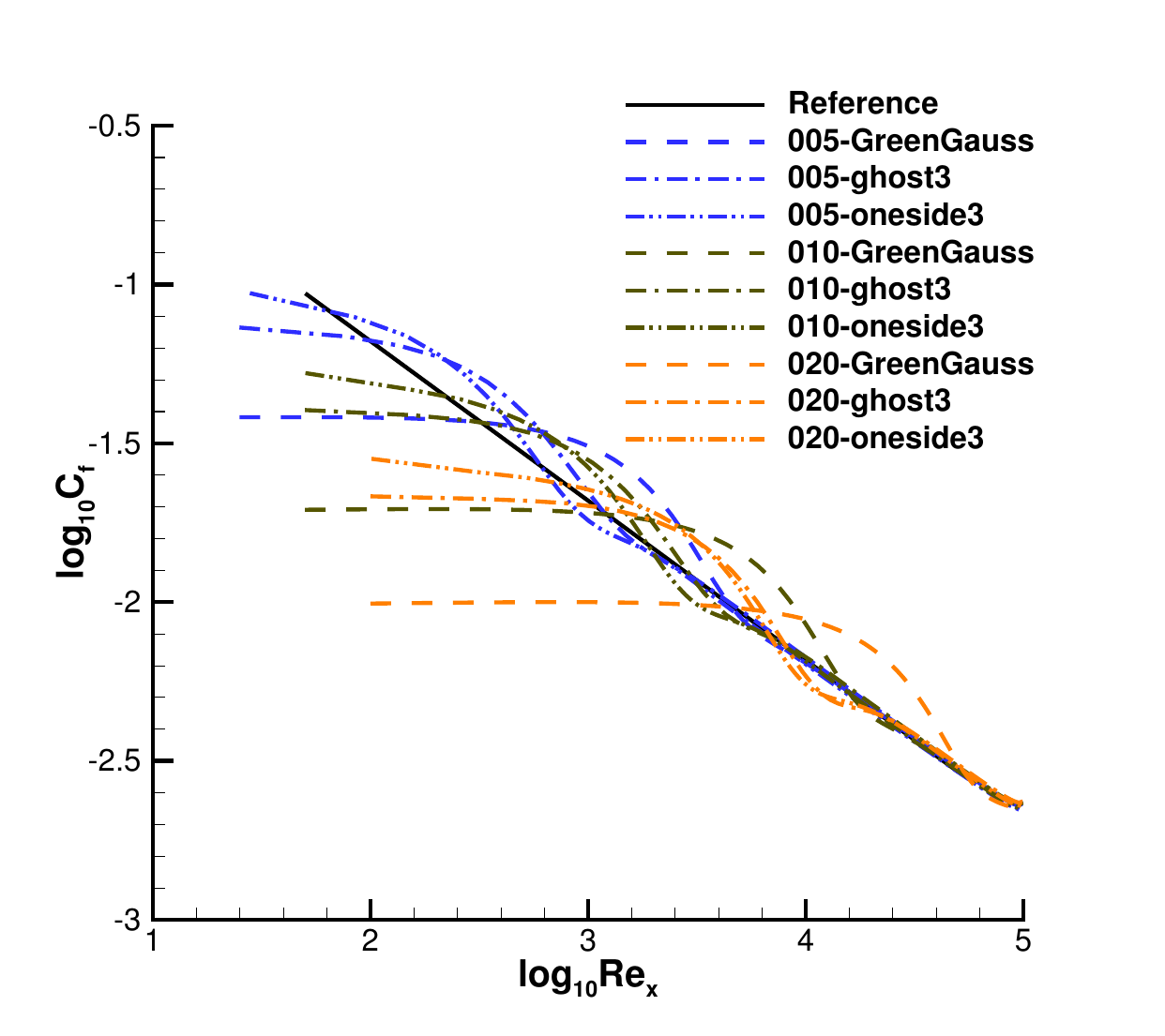}
		\caption{Distribution of skin-friction coefficient}
		\label{fig:bnLayer-cf}
	\end{figure}

	The distribution of skin-fraction coefficients $C_f$ along the plate is shown in Figure \ref{fig:bnLayer-cf}. Method \Rmnum{2} and Method \Rmnum{3} give much better estimation of the distribution of $C_f$ than Method \Rmnum{1}, and Method \Rmnum{3} performs better than Method \Rmnum{2}, especially in the front part of the plate.
	\subsection{Subsonic viscous flow around a cylinder at Re = 40}
\begin{figure}[!htbp]
		\centering
		\subfigure[Mesh]{\includegraphics[width = 0.45\columnwidth, trim = 0 50 0 50, clip]{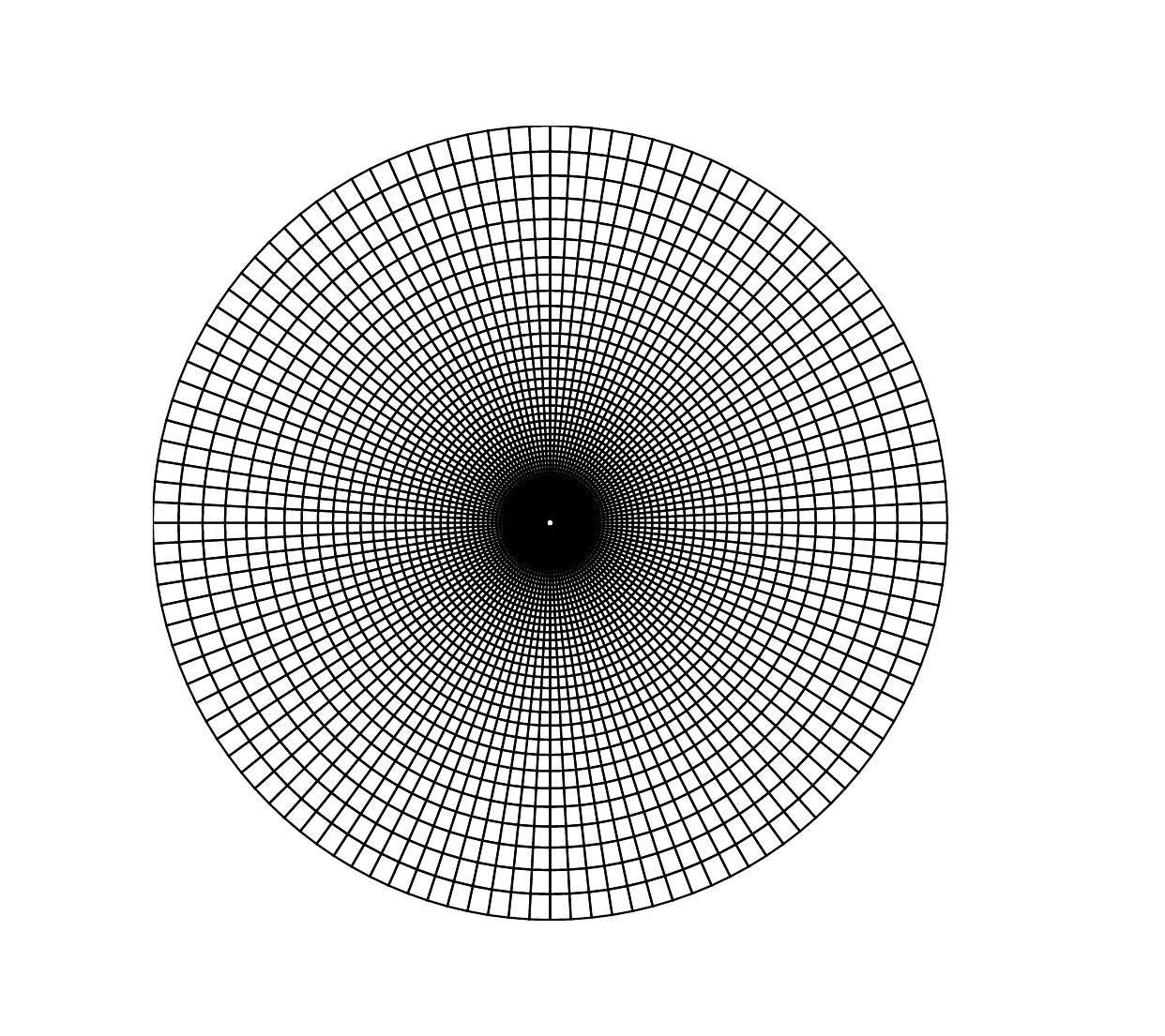}} \hspace{5pt}
		\subfigure[Streamline]{\includegraphics[width = 0.45\columnwidth, trim = 0 40 0 50, clip]{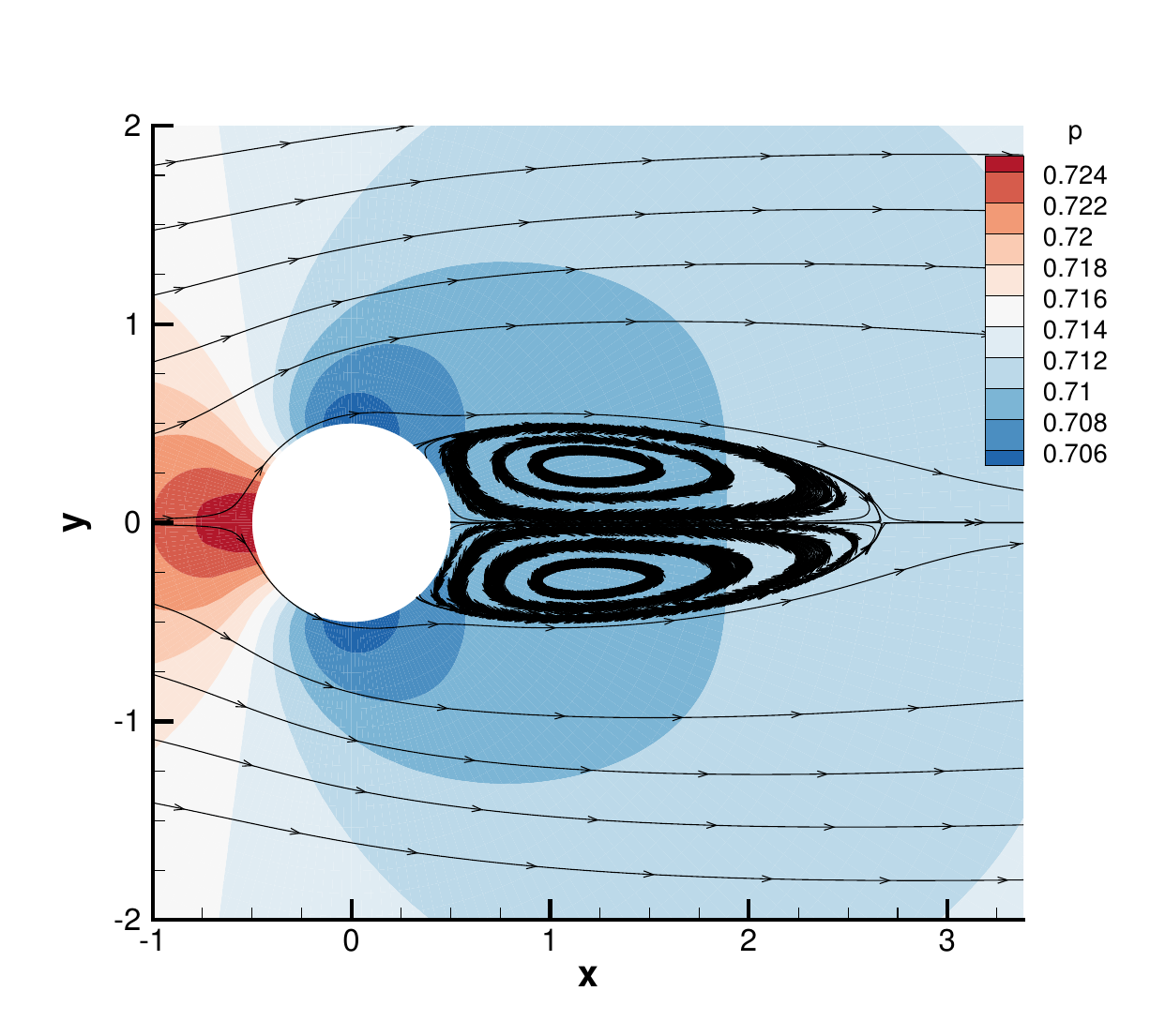}} \\
		\caption{Schematic of subsonic flow over a cylinder on $241\times 114$ mesh with $h = 1/96$ by Method \Rmnum{3}}
		\label{fig:cylinder-sch}
	\end{figure}

	A subsonic viscous flow around a cylinder is simulated. The incoming Mach number $\rm{Ma} = 0.15$, Reynolds number $\rm{Re} = 40$ and the characteristic length $L = 1$. The computational domain is shown in Figure \ref{fig:cylinder-sch}, where the diameter of the cylindrical computation domain is $D_{mesh} = 96.0$. Eighteen sets of grid are adopted to discretize the computation domain, as shown in Table \ref{tab:cylinder-set}.
	
	The non-slip adiabatic wall is imposed on the surface of the cylinder. The non-reflecting boundary condition based on the Riemann invariants is adopted on the outer boundary of the domain, and the incoming flow state is set as $\rho_{\infty} = 1, p_{\infty} = 1/\gamma$.
	The same numerical method as the laminar boundary layer case is used.
	
	The drag coefficient 
$$C_D = \frac{2 F_{drag}}{\rho_{\infty} U_{\infty}^2 L}$$ 
and wake length $l$ are listed in Table \ref{tab:cylinder-set}. The wake can not be resolved by Method \Rmnum{1} on the grids with $h = 1/6$, but still can be resolved by Method \Rmnum{2} and Method \Rmnum{3} on the same grids. And Method \Rmnum{3} performs better on coarse grids than Method \Rmnum{2}.
	\begin{table}[!htbp]
		\small
		\begin{center}
			\def\temptablewidth{1.0\textwidth}
			{\rule{\temptablewidth}{1pt}}
			\begin{tabular*}{\temptablewidth}{@{\extracolsep{\fill}}c|c|c|c|c|c|c|c}
				\multirow{2}{*}{Total mesh number} & \multirow{2}{*}{Near wall size} & \multicolumn{2}{l}{Method \Rmnum{1}} & \multicolumn{2}{l}{Method \Rmnum{2}} & \multicolumn{2}{l}{Method \Rmnum{3}} \\ \cline{3-8}
				&                    & $C_d$         & $L$         & $C_d$         & $L$         & $C_d$         & $L$         \\ \hline
				$241\times 114$ &1/96 & 1.526  & 2.26 & 1.526 & 2.25 &1.526 & 2.25 \\ 	
				$241\times 114$ &1/48 & 1.529  & 2.27 & 1.526 & 2.25 &1.526 & 2.25 \\
				$241\times 114$ &1/24 & 1.539  & 2.32 & 1.526 & 2.25 &1.526 & 2.23 \\
				$241\times 114$ &1/12 & 1.475  & 1.61 & 1.522 & 2.19 &1.529 & 2.18 \\
				$241\times 114$ &1/6 & 1.16  & -- & 1.463 & 1.63 &1.544 & 2.03 \\
				\hline
				$121\times 57$ &1/96 & 1.526  & 2.13 & 1.525 & 2.13 &1.525 & 2.13 \\ 	
				$121\times 57$ &1/48 & 1.529  & 2.15 & 1.526 & 2.13 &1.525 & 2.13 \\
				$121\times 57$ &1/24 & 1.541  & 2.21 & 1.526 & 2.13 &1.526 & 2.11 \\
				$121\times 57$ &1/12 & 1.477  & 1.54 & 1.523 & 2.07 &1.530 & 2.05 \\
				$121\times 57$ &1/6 & 1.151  & -- & 1.463 & 1.53 &1.544 & 1.91 \\
				\hline
				$61\times 29$ &1/96 & 1.510  & 1.47 & 1.508 & 1.47 &1.509 & 1.47 \\ 	
				$61\times 29$ &1/48 & 1.524  & 1.52 & 1.519 & 1.50 &1.520 & 1.50 \\
				$61\times 29$ &1/24 & 1.546  & 1.59 & 1.527 & 1.52 &1.528 & 1.51 \\
				$61\times 29$ &1/12 & 1.485  & 1.15 & 1.530 & 1.48 &1.537 & 1.46 \\
				$61\times 29$ &1/6 & 1.133  & -- & 1.468 & 1.12 &1.551 & 1.36 \\
				\hline
				$33\times 33$ &1/24 & 1.642  & 1.03 & 1.622 & 0.99 &1.625 & 0.99 \\
				$33\times 33$ &1/12 & 1.567  & 0.72 & 1.622 & 0.97 &1.633 & 0.97 \\
				$33\times 33$ &1/6 & 1.236  & -- & 1.559 & 0.78 &1.642 & 0.92 \\
			\end{tabular*}
			{\rule{\temptablewidth}{1pt}}
		\end{center}
		\vspace{-4mm} \caption{\label{viscous subsonic sphere} Quantitative results of the cylinder at Re = 40}
		\label{tab:cylinder-set}
	\end{table}
	
\begin{figure}[!htbp]
		\centering
		\subfigure{\includegraphics[width = 0.45\columnwidth, trim = 0 10 0 10, clip]{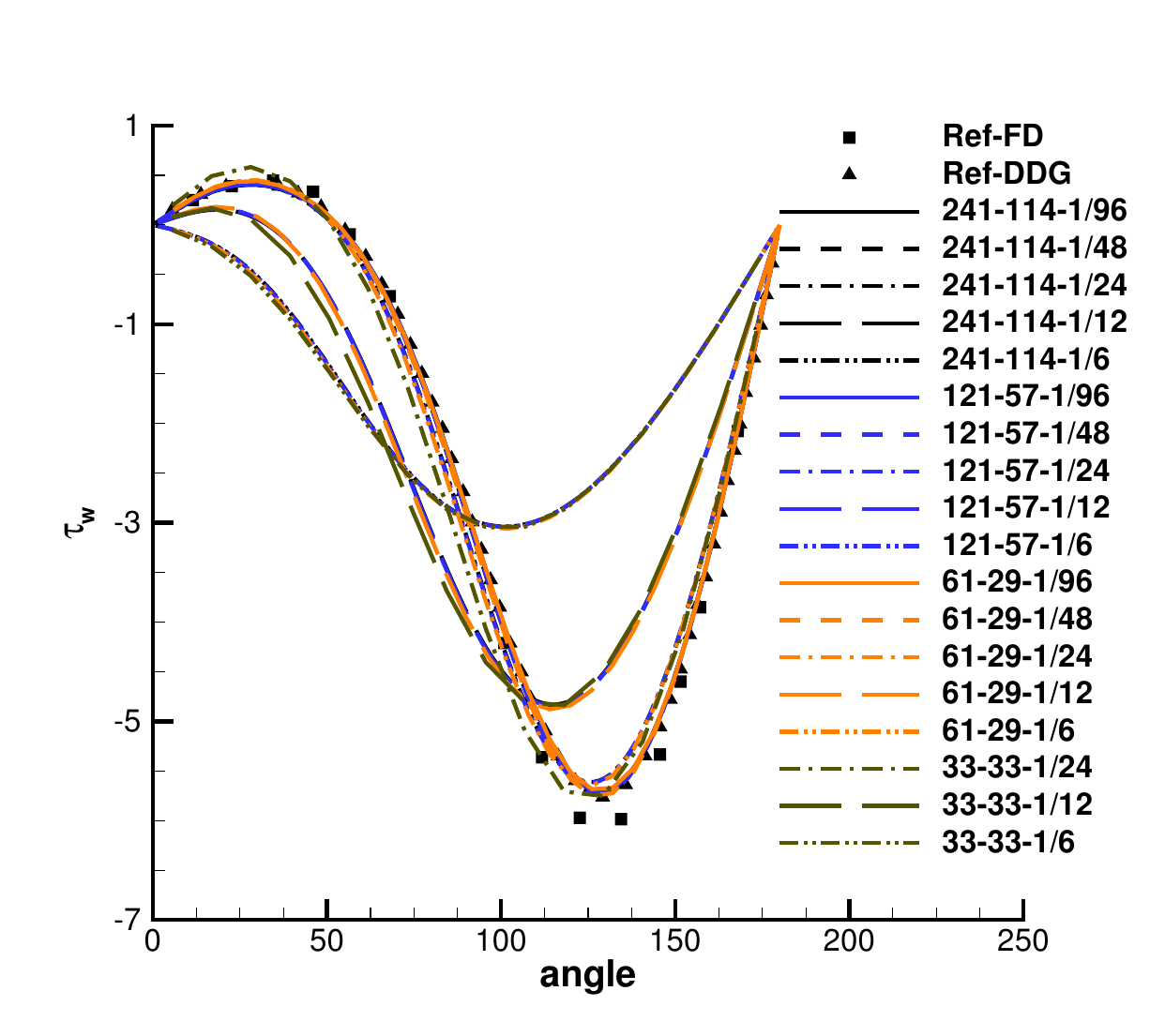}} \hspace{5pt}
		\subfigure{\includegraphics[width = 0.45\columnwidth, trim = 0 10 0 10, clip]{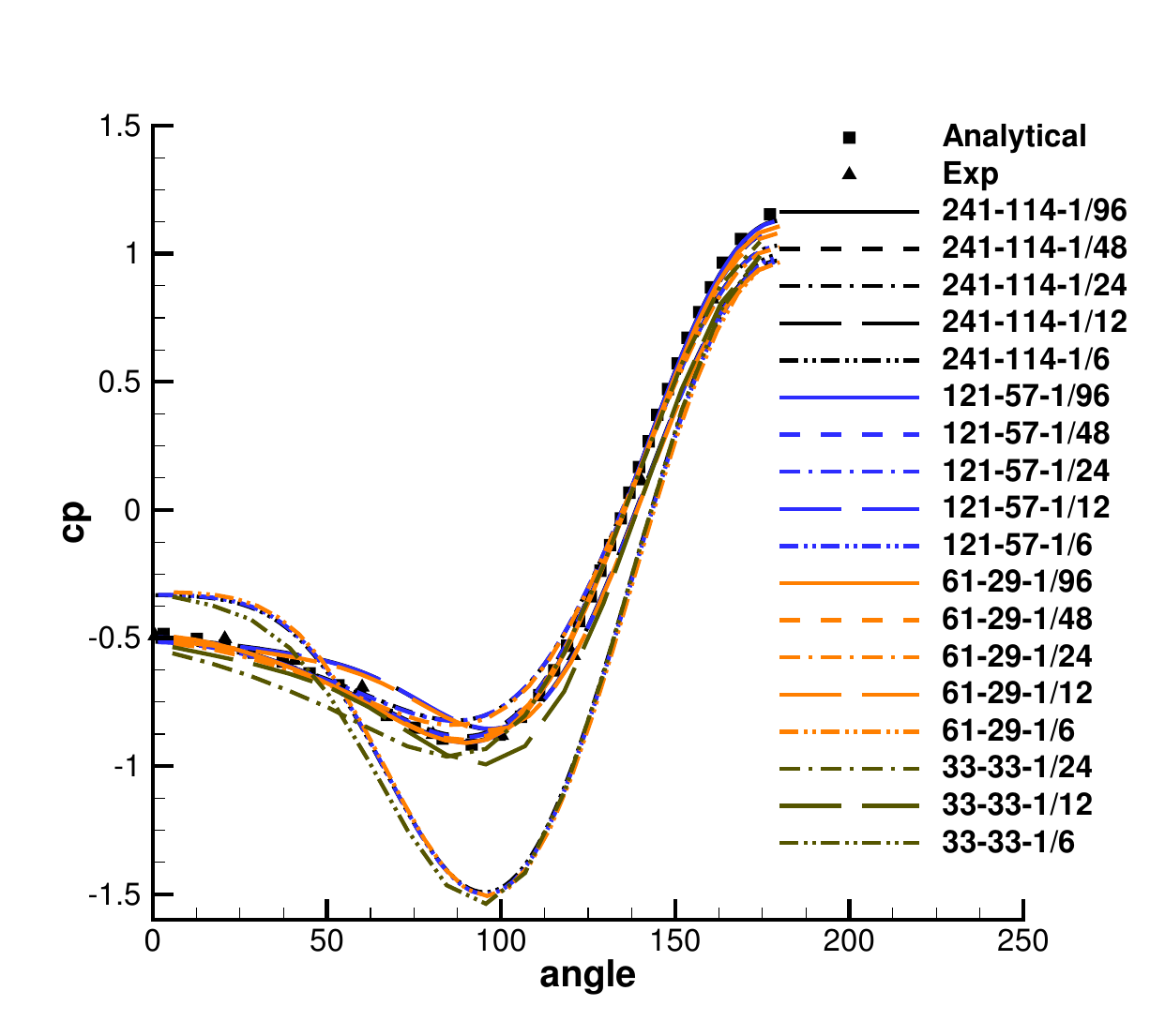}} \\
		\subfigure{\includegraphics[width = 0.45\columnwidth, trim = 0 10 0 10, clip]{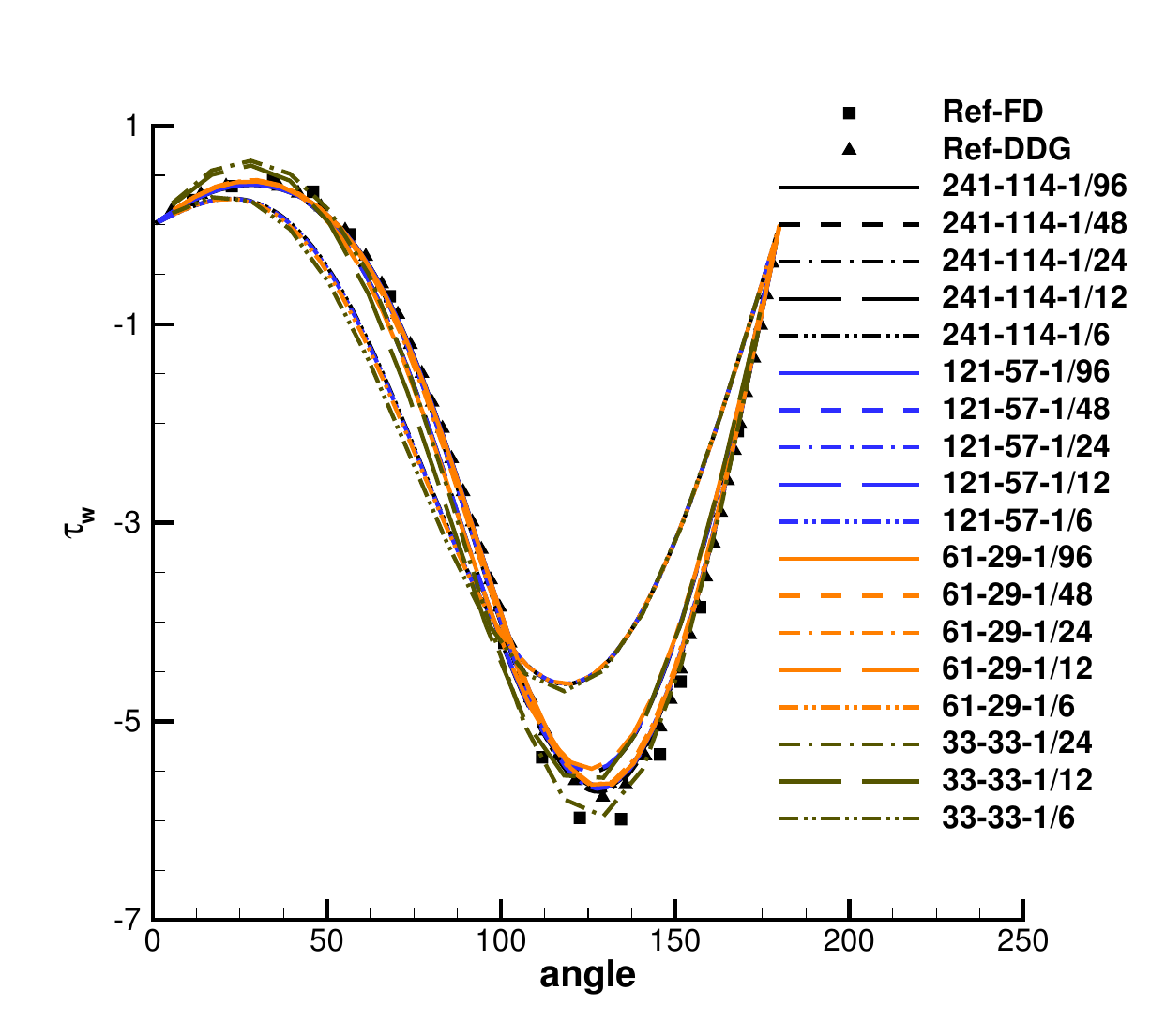}} \hspace{5pt}
		\subfigure{\includegraphics[width = 0.45\columnwidth, trim = 0 10 0 10, clip]{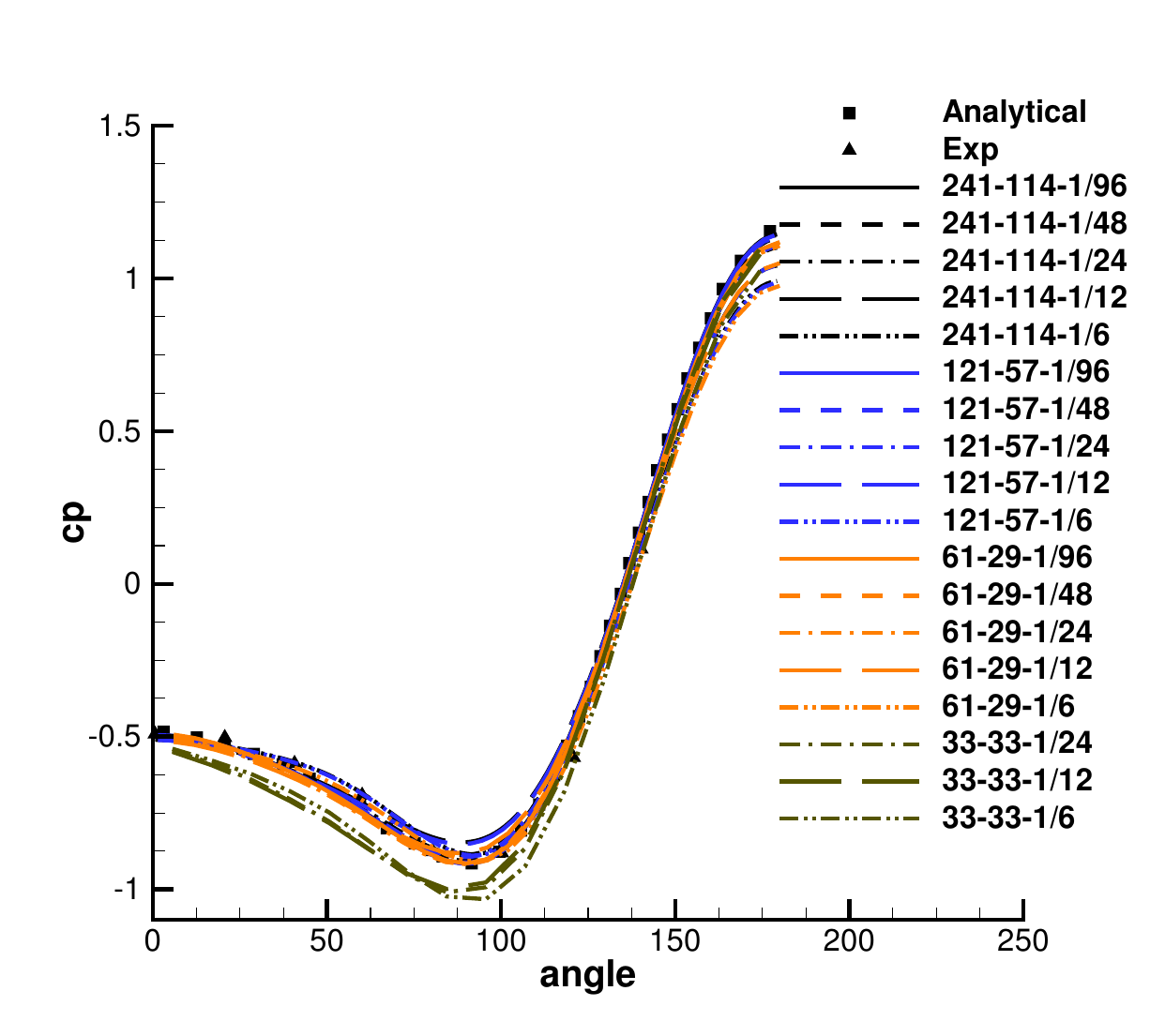}} \\
		\subfigure{\includegraphics[width = 0.45\columnwidth, trim = 0 10 0 10, clip]{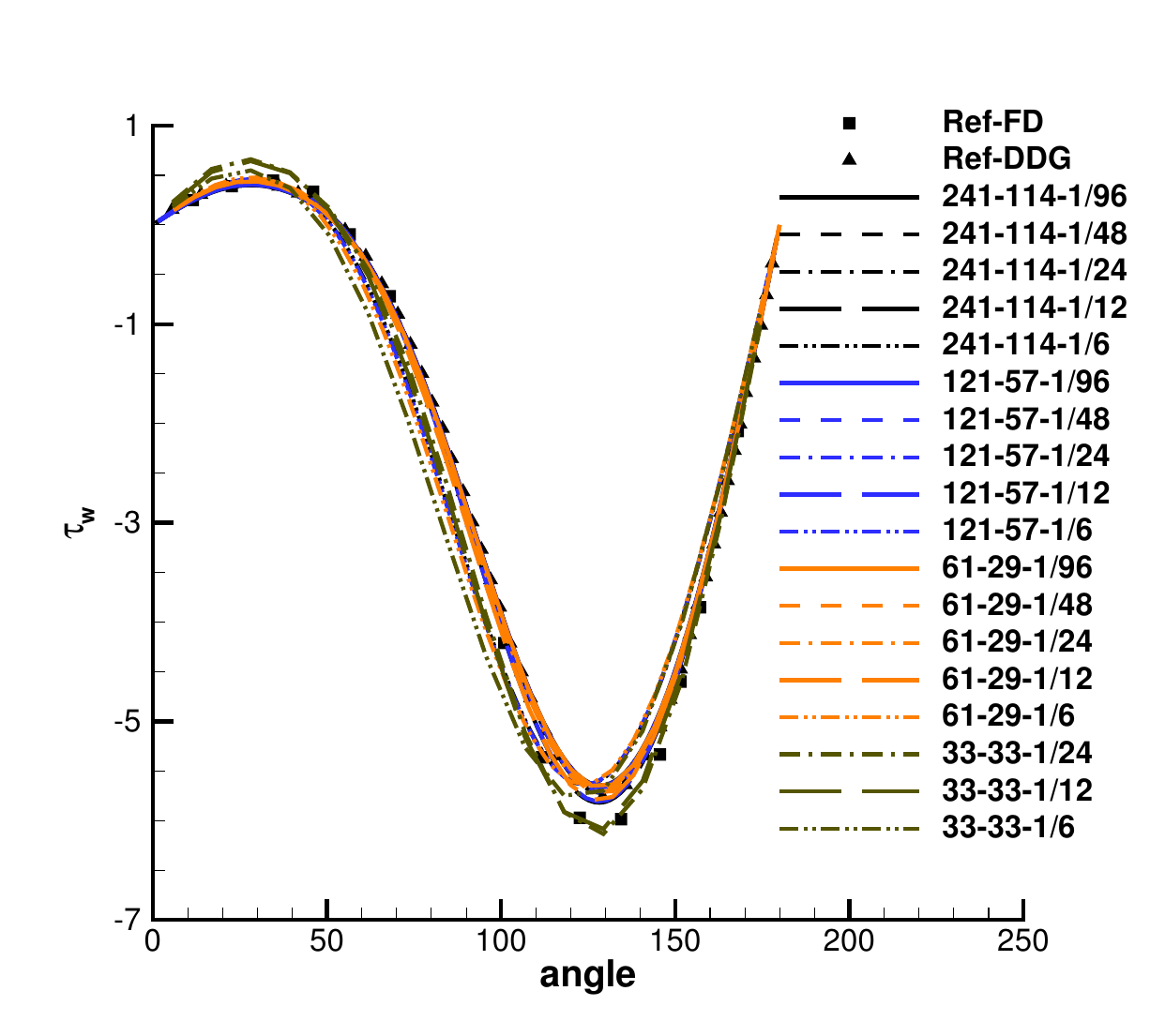}} \hspace{5pt}
		\subfigure{\includegraphics[width = 0.45\columnwidth, trim = 0 10 0 10, clip]{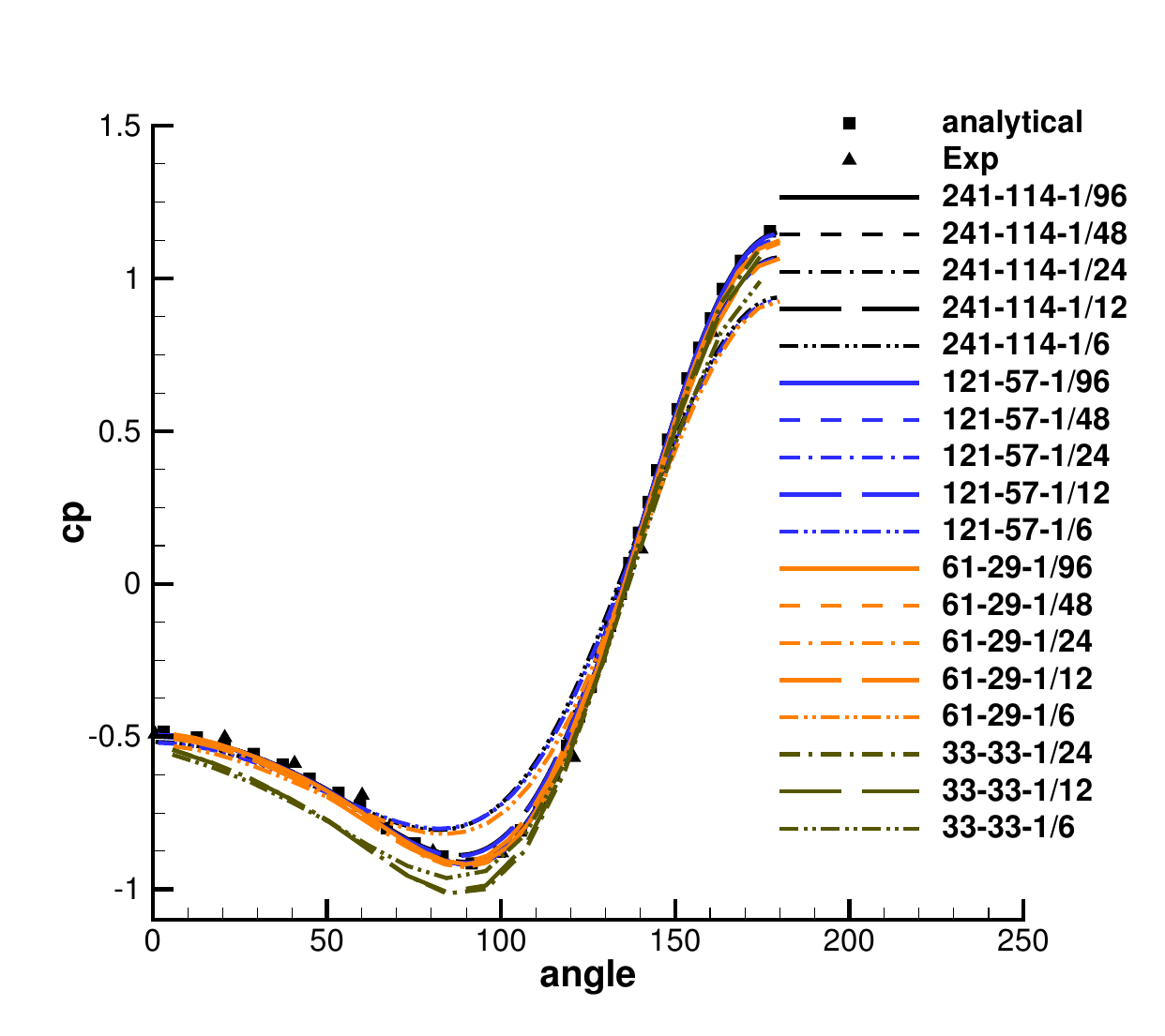}} \\
		\caption{Cylinder: Re = 40. Left: surface local tangential velocity gradient distribution, right: surface pressure coefficient distribution. Top row: Method \Rmnum{1}, middle row: Method \Rmnum{2}, bottom row: Method \Rmnum{3}}
		\label{fig:cylinder-tauw}
	\end{figure}

	The distributions of pressure coefficient 
$$C_p=\frac{2(p-p_{\infty})}{\rho_{\infty} U_{\infty}^2}$$ 
and local tangential velocity gradient 
$$\tau_w = \frac{D}{2U_{\infty} }\frac{\partial U_{\tau}}{\partial \eta}$$ 
along the cylinder surface are shown in Figure \ref{fig:cylinder-tauw}. Similar with the laminar boundary layer case, the third-order Method \Rmnum{2} and Method \Rmnum{3} can achieve much better results than the second-order Method \Rmnum{1}. The results of Method \Rmnum{1} deviate the reference severely on grids with $h = 1/12$, but Method \Rmnum{2} still obtains a reasonable result on grids with $h = 1/12$ but fails on grids with $h = 1/6$, and the results of Method \Rmnum{2} on grids with $h = 1/6$ is very close to those of Method \Rmnum{1} on grids with $h = 1/12$, and Method \Rmnum{3} obtains a rather good result even on grids with $h = 1/6$. It is shown in Figure \ref{fig:cylinder-tauw} that the quality of the solution at wall boundary is closely related to the height of the first grid layer, where similar results will be obtained if the same first layer grid height is used, as shown in Table \ref{tab:cylinder-set}.  
The length of the wake depends mainly on the grid size of interior cells in the computational domain. When the boundary layer is not well resolved, i.e., on the grids with $h = 1/6$, Method \Rmnum{3} predicts more accurate wake length than Method \Rmnum{2}.

	\subsection{Viscous hypersonic flow around a cylinder}
\begin{figure}[!htbp]
		\centering
		\subfigure[Mesh]{\includegraphics[width = 0.15\columnwidth]{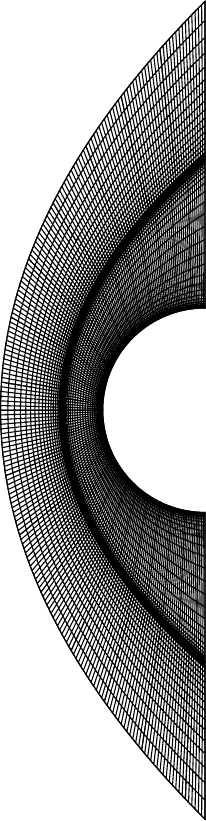}} \hspace{30pt}
		\subfigure[Density]{\includegraphics[width = 0.15\columnwidth]{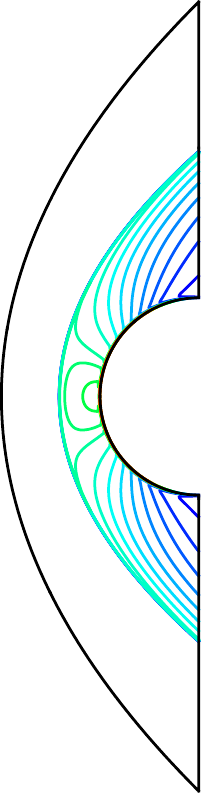}} \hspace{30pt}
		\subfigure[Pressure]{\includegraphics[width = 0.15\columnwidth]{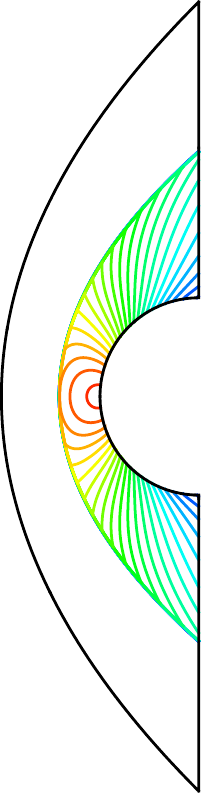}} \hspace{30pt}
		\subfigure[Mach]{\includegraphics[width = 0.15\columnwidth]{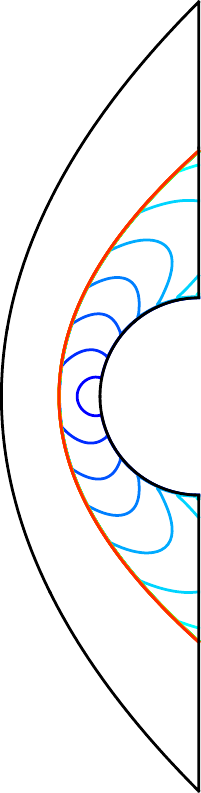}} \\
		\caption{Viscous hypersonic flow over cylinder with Ma = 8.03.(a) Mesh, (b) Density, (c) Pressure, (d) Mach number}
		\label{fig:halfCylinder}
	\end{figure}

	In this section, viscous hypersonic flows around a cylinder with non-slip isothermal boundary are tested. The incoming Mach number $\rm{Ma} = 8.03$, Reynolds number $\rm{Re} = 1.835 \times 10^5$, Prandtl number $Pr = 0.72$, $T_{\infty} = 124.94K$ and the characteristic length $L = 1$. The computational domain is shown in Figure \ref{fig:halfCylinder}(a),  the diameter of the inner cylindrical wall is $D = 1$. The computation domain is discretized  into $400\times 110$  quadrilateral cells with the height of the first layer grid $h= 10^{-4}$. For a clear presentation, the mesh shown in Figure \ref{fig:halfCylinder}(a) has $160\times 80$ mesh points. 
To capture shock and resolve boundary layer better, meshes near shock front and boundary wall are fined locally.
	
	The non-slip isothermal wall is imposed on the surface of the cylinder, where the wall temperature is $T_w = 294.44K$. The non-reflecting boundary condition based on the Riemann invariants is adopted on the outer boundary of the domain, and the outlet is set as supersonic outflow. The nonlinear reconstruction Eq. \ref{nonlinear-p2} and the complete gas-kinetic solver Eq. \ref{f_inter} are applied in this case, and the gas-kinetic solver for isothermal boundary introduced in Section \ref{gks_iso} is used to evaluate the flux across isothermal wall. CFL number is taken as 0.1 in this case.
	
\begin{figure}[!htbp]
		\centering
		\subfigure[Pressure]{\includegraphics[width = 0.45\columnwidth]{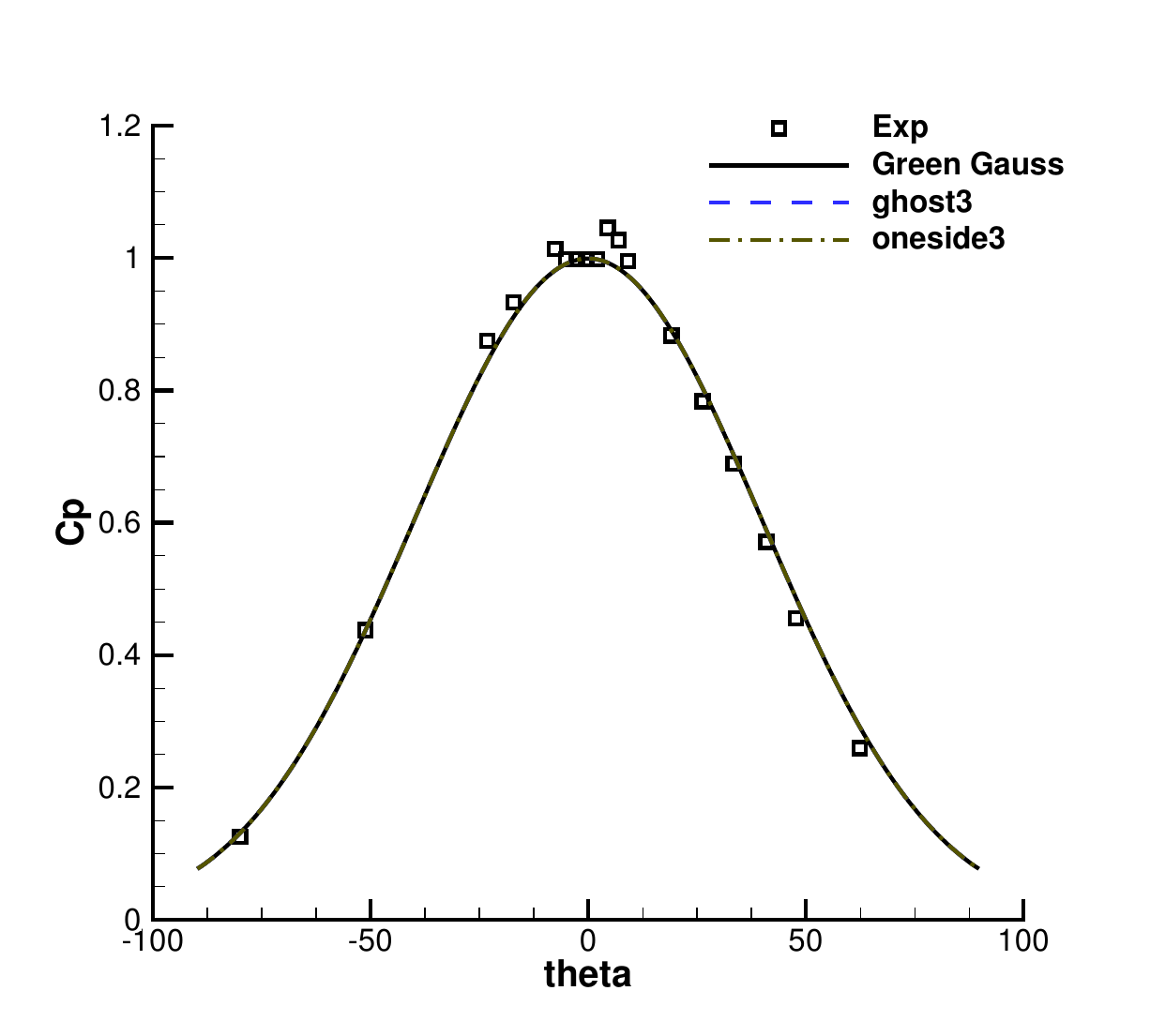}} \hspace{5pt}
		\subfigure[Heat flux]{\includegraphics[width = 0.45\columnwidth]{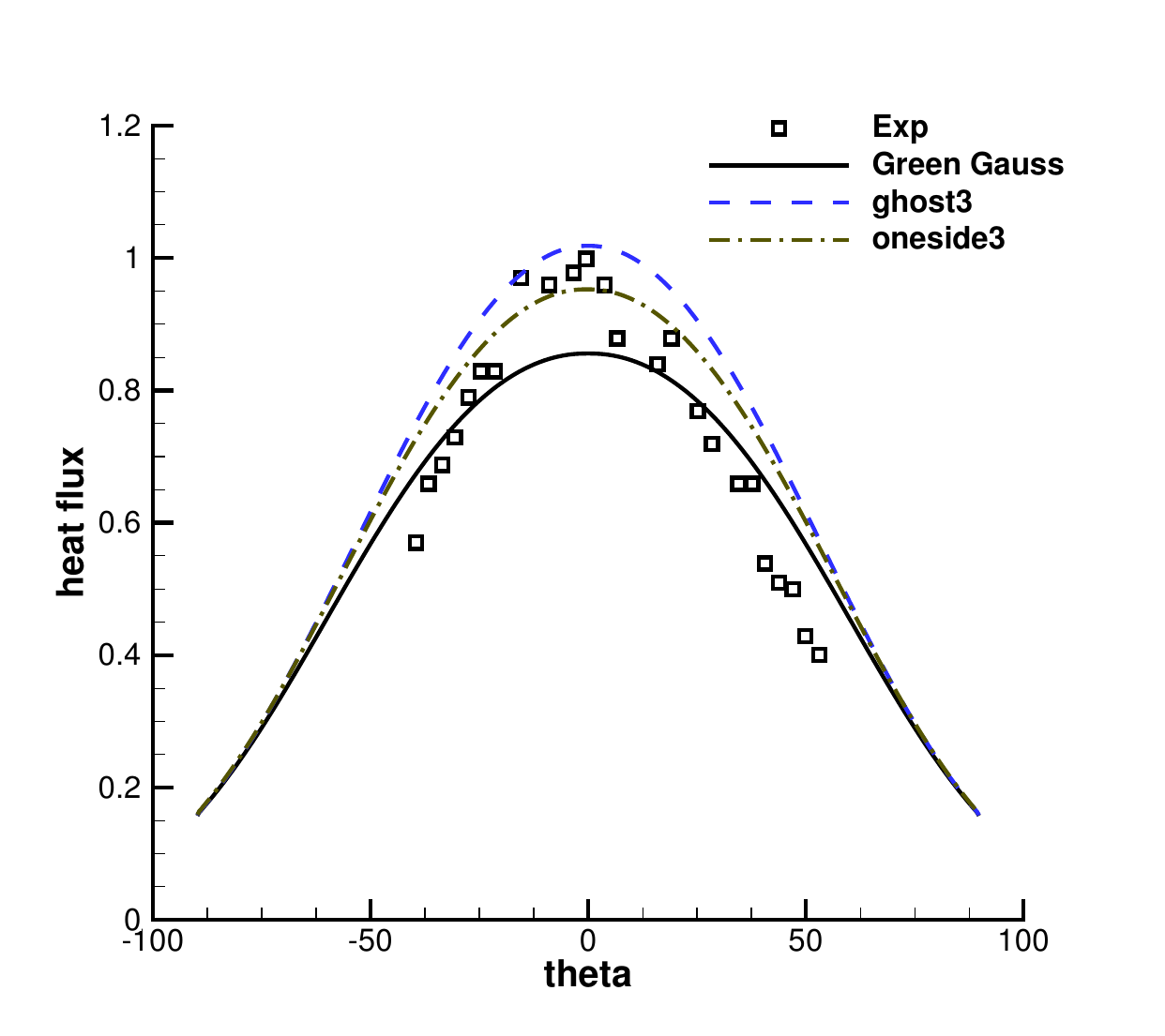}}
		\caption{The distribution of normalized pressure and heat flux along the cylinder}
		\label{fig:ma8-surface}
	\end{figure}

	Figure \ref{fig:halfCylinder}(b-d) demonstrate the contours of density, pressure and Mach number. The pressure and heat flux are normalized by $0.9209\rho_{\infty} U_{\infty}^2$ and $0.003655\rho_{\infty} U_{\infty}^3$ respectively. The normalized pressure and heat flux along the surface of the cylinder are compared with experimental data in Ref. \cite{cylinder-exp} in Figure \ref{fig:ma8-surface}. The pressure on the surface of the cylinder can be captured accurately, the one-side third-order boundary treatment gives the most accurate prediction of heat flux, with a large first layer mesh size with cell Reynolds number $\approx 18.35$. 
This case validates the effectiveness of the high-order boundary treatment for hypersonic thermodynamical problems.
	
\begin{figure}[!htbp]
		\centering
		\includegraphics[width = 0.5\columnwidth]{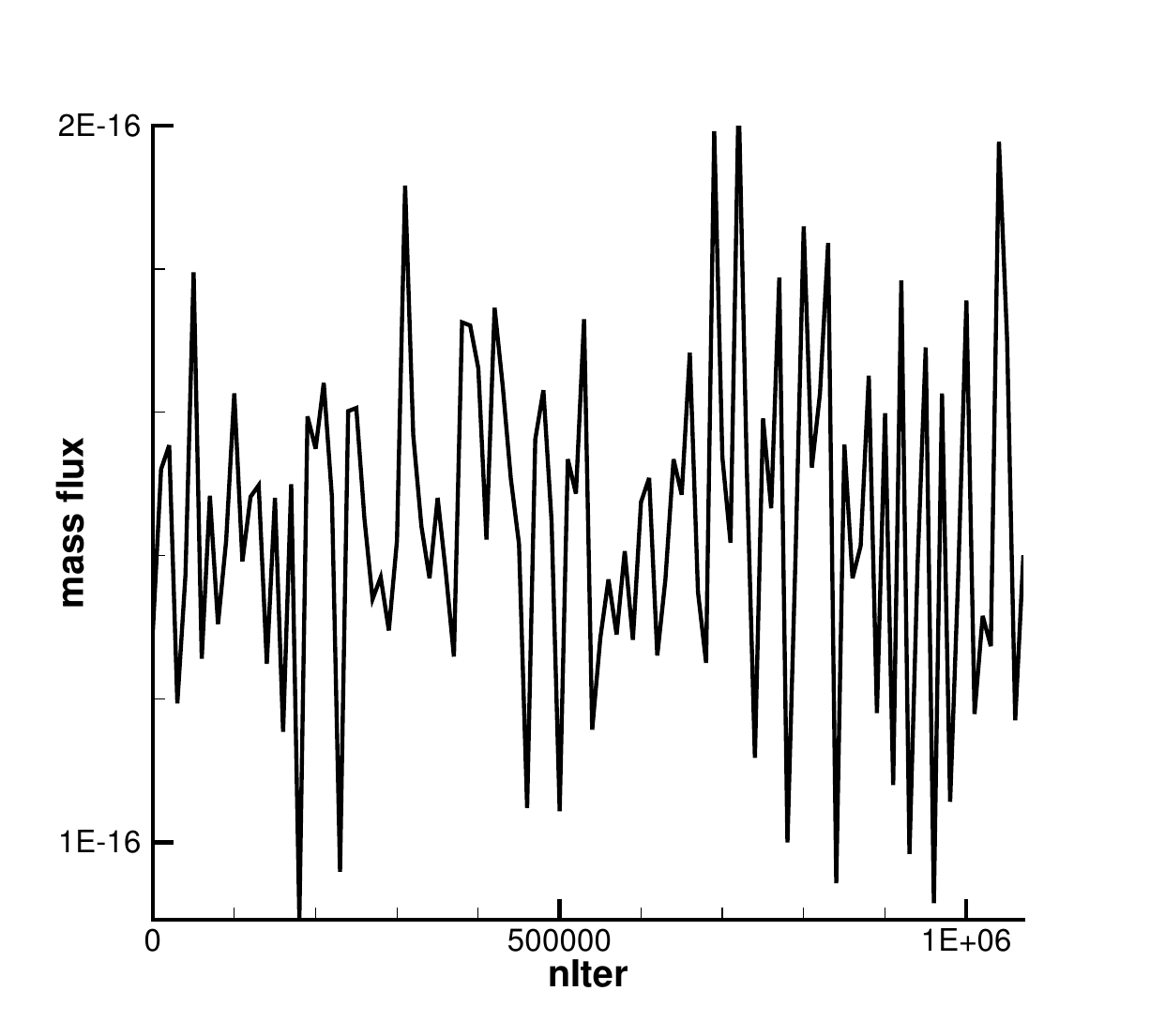}
		\caption{Mass flux across the cylinder surface during iteration}
		\label{fig:massflux}
	\end{figure}

	Figure \ref{fig:massflux} demonstrates the mass flux across the wall boundary during the iteration calculated by Method \Rmnum{3}, the mass flux is written every 10000 iterations. Due to the gas-kinetic isothermal flux introduced in Section \ref{gks_iso} and the reconstruction for boundary cells in Section \ref{recon_iso}, the mass flux keeps at the level of machine zero, even at the very beginning of the iteration. This property is vitally important for isothermal boundary, or the mass leakage through the wall boundary will lead to underestimation of heat flux.
	
	\section{Conclusion}
	In this paper, a class of one-sided third-order compact gas-kinetic scheme for non-slip wall boundaries is constructed. The method is validated through subsonic to hypersonic aerodynamic problems. The current boundary treatments are high-order, simple, and easy extension to curved boundaries.
	The updated flow variables and their gradients from compact GKS are used in the construction of the numerical non-slip boundary discretization with a third-order  spatial accuracy in the boundary cells. Based on the  time-accurate flux function in the GKS, two-step fourth-order temporal discretization is adopted to achieve high-order temporal accuracy.
	A kinetic boundary condition for isothermal wall is further  developed to guarantee the no-penetration condition through the wall boundary. The satisfaction of no mass penetration is vitally important for the accurate estimation of heat flux. In the numerical tests, the third-order compact scheme with ghost cells and the proposed one-sided third-order compact scheme on curved meshes demonstrate significantly better performance than the second-order boundary condition discretization through Green-Gauss method. Furthermore, the proposed one-sided third-order compact scheme on curved meshes yields superior results compared to the third-order scheme with ghost cells, particularly when the boundary layer is not well-resolved for the flow computation with curved boundaries.

	\section*{Acknowledgements}
	
	We are grateful to Dr Yajun Zhu for numerical discussions during this work and helpful comments on this manuscript. The current research is supported by the National Key R\&D Program of China 2022YFA1004500,
	National Natural Science Foundation of China (Grant No. 12172316, 92371107),
	and Hong Kong research grant council (16301222, 16208324).

	\bibliographystyle{ieeetr}
	\bibliography{reference}
	
	\begin{appendices}
		\section{Quadratic integration on curved line}
		The line-averaged quantities in Eq. \ref{line-aver} and \ref{line-grad-aver} is calculated in the reference system based on the isoparametric transformation. 
    The isoparametric transformation of the quadratic curved line Eq. \ref{line-iso} is
		\begin{equation}
			\dot{x} = \sum_{k=1}^{3} v_k^{'}(\xi) \mathbf{x}_k = (-3x_0-x_1+4x_2) + (4x_0+4x_1-8x_2)\xi = a_x + b_x \xi, \nonumber
		\end{equation}
		\begin{equation}
			\dot{y} = \sum_{k=1}^{3} v_k^{'}(\xi) \mathbf{y}_k = (-3y_0-y_1+4y_2) + (4y_0+4y_1-8y_2)\xi = a_y + b_y \xi. \nonumber
		\end{equation}
		The unit length of the curve is
		\begin{equation}
			\mathrm{d}s = \sqrt{\dot{x}^2 + \dot{y}^2}\mathrm{d}\xi = \sqrt{\phi_0(\xi + \phi_1 )^2 + \phi_2}\mathrm{d}\xi,
		\end{equation}
		where
		\begin{equation}
			\phi_0 = b_x^2 + b_y^2, \nonumber
		\end{equation}
		\begin{equation}
			\phi_1 = \frac{a_xb_x + a_yb_y}{ b_x^2 + b_y^2}, \nonumber
		\end{equation}
		\begin{equation}
			\phi_2 = a_x^2 + a_y^2 - \frac{(a_xb_x + a_yb_y)^2}{b_x^2 + b_y^2}. \nonumber
		\end{equation}
		
		To calculate the curve-averaged value of the basis function of the $p^2$ polynomial in Eq. \ref{p_2}, the following integral needs to be calculated
		\begin{equation}
			\int_{0}^{1} (\xi + \phi_1)^p \sqrt{\phi_0(\xi + \phi_1 )^2 + \phi_2}\rm{d}\xi = (\frac{\phi_2}{\phi_0})^{\frac{p+1}{2}} \sqrt{\phi_2} \int_{\theta_0}^{\theta_1} tan^p \theta sec^3 \theta \rm{d} \theta,
		\end{equation}
		where $\theta_0 = \arctan(\sqrt{\frac{\phi_0}{\phi_2}} \phi_1)$, $\theta_1 = \arctan(\sqrt{\frac{\phi_0}{\phi_2}} (1 + \phi_1))$. If $p$ is an odd integer $p=2k+1$,
		\begin{equation}
			\int \tan^{2k+1} \theta \sec^3 \theta \rm{d} \theta = \int (sec^2 \theta - 1)^k sec^2 \theta \rm{d} \sec \theta.
		\end{equation}
		If $p$ is an even integer $p=2k$,
		\begin{equation}
			\int \tan^{2k} \theta \sec^3 \theta \rm{d} \theta = \int (sec^2 \theta - 1)^k sec^2 \theta \rm{d}  \theta,
		\end{equation}
		and
		\begin{equation}
			I_n = \int \sec^n \theta \rm{d} \theta = \frac{1}{n - 1} (tan\theta \sec^{n - 2} \theta + (n - 2)I_{n-2}). \nonumber
		\end{equation}
		
		\section{Quadratic integral on curved triangle}
		The isoparametric transformations of the quadratic triangular cell Eq. \ref{tri-iso} is
		\begin{equation}
			\begin{aligned}
				\frac{\partial x}{\partial \xi} = a_0\xi + a_1\eta + a_2,  \\
				\frac{\partial x}{\partial \eta} = a_3\xi + a_4\eta + a_5, \\
				\frac{\partial y}{\partial \xi} = b_0\xi + b_1\eta + b_2,  \\
				\frac{\partial y}{\partial \eta} = b_3\xi + b_4\eta + b_5, \\
			\end{aligned}
		\end{equation}
		where $a_0 = 4(x_0 + x_1 - 2x_3), a_1 = 4(x_0 - x_3 + x_4 - x_5), a_2 = - 3x_0 - x_1 + 4x_3$,$a_3 = 4(x_0 - x_3 + x_4 - x_5), a_4 = 4(x_0 + x_2 - 2x_5), a_5 = -3x_0 - x_2 + 4x_5$,$b_0 = 4(y_0 + y_1 - 2y_3), b_1 = 4(y_0 - y_3 + y_4 - y_5), b_2 = - 3y_0 - y_1 + 4y_3$,$b_3 = 4(y_0 - y_3 + y_4 - y_5), b_4 = 4(y_0 + y_2 - 2y_5), b_5 = -3y_0 - y_2 + 4y_5$.\\
		The determinant of the Jocabian matrix of the transform is
		\begin{equation}
			\begin{aligned}
				\left| \begin{array} {cc}
					\frac{\partial x}{\partial \xi} & \frac{\partial x}{\partial \eta} \\
					\frac{\partial y}{\partial \xi} & \frac{\partial y}{\partial \eta}
				\end{array}
				\right|
				&=a_2b_5 - a_5b_2 + (a_0b_5 + a_2b_3 - a_3b_2 - a_5b_0)\xi \\
				&+ (a_1b_5 + a_2b_4 - a_4b_2 - a_5b_1)\eta + (a_0b_3 - a_3b_0)\xi^2 \\
				&+ (a_0b_4 + a_1b_3 - a_3b_1 - a_4b_0)\xi \eta + (a_1b_4 - a_4b_1) \eta^2.
			\end{aligned} \nonumber
		\end{equation}
		The integrals of the basis functions on the curved triangle in reference system are
		\begin{equation}
			\int_0^1 \int_{0}^{1-\xi} \xi^{m} \eta^{n} \mathrm{d}\eta \mathrm{d} \xi = \frac{m!n!}{(m+n+2)!}.
		\end{equation}
	\end{appendices}
	
\end{document}